%Paper: alg-geom/9406008
%From: cornalba@math.ias.edu (Maurizio D. Cornalba)
%Date: Thu, 30 Jun 94 13:18:02 EDT

%The text following the \enddocument command is a uuencoded,
%tar-compressed postscript file containing the figures. To print them,
%strip off the rest of the message, uudecode, decompress, and send
%resulting psfile to printer.

\def\mgn#1#2{{\Cal M}_{#1,#2}}
\def\mgnbar#1#2{\overline{\Cal M}_{#1,#2}}
\def\mcomb#1#2{{\Cal M}_{#1,#2}^{comb}}
\def\mcombar#1#2{\overline{\Cal M}_{#1,#2}^{comb}}
\def\mprime#1#2{\overline{\Cal M}_{#1,#2}'}
\def\mm{\overline{\Cal M}}
\def\mp{\overline{\Cal M}'}
\def\wmn#1#2{W_{#1,#2}}
\def\ul#1{\underline{#1}}
\def\trace{\operatorname{tr}}
\def\ang#1{\langle#1\rangle}
\def\br#1{\langle#1\rangle_\Lambda}

\def\brr#1{\langle\langle#1\rangle\rangle_\Lambda}
\def\bbrr#1{\left\langle\left\langle#1\right\rangle\right\rangle_\Lambda}

\def\data{6/16/1994}

\font\titolo=cmbx12

\documentstyle{amsppt}
\magnification 1200
\pagewidth{16.5truecm}
\pageheight{22.5truecm}
\NoBlackBoxes
\leftheadtext\nofrills{\it E. Arbarello, M. Cornalba, Cohomology classes on
moduli --
version \data}
\rightheadtext\nofrills{\it E. Arbarello, M. Cornalba, Cohomology classes on
moduli
-- version \data}

\document
\centerline{\titolo Combinatorial and algebro-geometric
cohomology classes}
\centerline{\titolo on the moduli spaces of curves}
\centerline{(version \data)}
\vskip0.8truecm
\centerline{Enrico Arbarello\footnote{\sevenrm Supported in
part by the M.U.R.S.T. 40\% project ``Geometria Algebrica''.}}
\vskip0.15truecm
\centerline{\it Dipartimento di Matematica Istituto ``Guido
Castelnuovo''}
\centerline{\it Universit\`a di Roma ``La Sapienza''}
\centerline{\it P.le A. Moro 2, 00185 Roma, Italia}
\centerline{\it and Institute for Advanced Study}
\centerline{\it Princeton, New Jersey 08540, U.S.A.}
\centerline{\it e-mail: arbarello\@ sci.uniroma1.it}
\bigskip
\centerline{Maurizio Cornalba$^{1,}$\footnote{\sevenrm Supported
in part by a grant from the Ellentuck Fund.}}
\vskip0.15truecm
\centerline{\it Dipartimento di Matematica, Universit\`a di Pavia}
\centerline{\it via Abbiategrasso 209, 27100 Pavia, Italia}
\centerline{\it and Institute for Advanced Study}
\centerline{\it Princeton, New Jersey 08540, U.S.A.}
\centerline{\it e-mail: cornalba\@ ipv85.unipv.it}

\vskip1.5truecm
\specialhead Introduction
\endspecialhead
\noindent Given a compact Riemann surface $C$, $n$ points on it,
and $n$ positive real numbers ($n>\max(0,2-2g)$), Strebel's theory of quadratic
differentials [15] provides a canonical way of dissecting $C$
into $n$ poligons and assigning lengths to their sides. As Mumford
first noticed, this can be used to give a combinatorial description of
the moduli space $\mgn{g}{n}$ of
$n$-pointed smooth curves of given genus $g$. If one looks at
moduli spaces from this point of view, one can construct
combinatorial cycles in them (cf. [7], for instance). It is then
natural to ask how these may be related to the algebraic geometry
of moduli space. It was first conjectured by Witten that the
combinatorial cycles can be expressed in terms of
Mumford-Morita-Miller classes. The first result in this direction is
due to Penner [13]; we will comment on his work at the
end of section 2. As we shall briefly explain now, and more
extensively in section 3, our approach to the question has its origin
in the papers [16], [17] and [7] by Witten and
Kontsevich. The combinatorial cycles we are talking about will be
denoted by the symbols
$\wmn{m_*}{n}$, where
$m_*=(m_0,m_1,m_2,\dots)$ is an infinite sequence of non-negative
integers, almost all zero, and $n$ a positive integer. On the moduli
space
$\mgn{g}{n}$ live particular cohomology classes of degree two,
denoted
$\psi_i$, $i=1,\dots,n$; by definition, $\psi_i$ is the Chern class of
the line bundle whose fiber at the point $[C; x_1,\dots,x_n]\in
\mgn{g}{n}$ is the cotangent space to $C$ at $x_i$. For the
intersection numbers of the $\psi_i$ along the $\wmn{m_*}{n}$ one
uses the notation
$$
\ang{\tau_{d_1}\dots\tau_{d_n}}_{m_*}=
\int_{\wmn{m_*}{n}}\prod_{i=1}^n\psi_i^{d_i}\,.
$$
These numbers are best organized as coefficients of an infinite
series
$$
F(t_0,t_1,\dots,s_0,s_1,\dots)=\sum\ang{\tau_{d_1}\dots\tau_{d_n}}_{m_*}
t_{d_1}\dots t_{d_n}\prod s_i^{m_i}\,.
$$
Kontsevich proves that $\exp F$ is an asymptotic expansion, as
$\Lambda^{-1}$ goes to zero, of the integral
$$
{\int_{\Cal{H}_N}\exp\left(-\sqrt{-1}\sum_{j=0}^\infty\left(-{1\over
2}\right)^j s_j {\trace(X^{2j+1})\over 2j+1}\right)\exp\left(-{1\over
2}\operatorname{tr}(X^2\Lambda)\right)dX
\over
\int_{\Cal{H}_N}\exp\left(-{1\over
2}\operatorname{tr}(X^2\Lambda)\right)dX}
\,,
$$
where $\Cal{H}_N$ is the space of $N\times N$ hermitian
matrices, $dX$ is a $U(N)$-invariant measure, and $\Lambda$ is a
positive definite diagonal $N\times N$ matrix, linked to the $t$
variables by the substitution $t_i=-(2i-1)!!\trace(\Lambda^{-2i-1})$.
Using this result, Di Francesco, Itzykson and Zuber [3]
showed that the derivatives of $\exp F$ with respect to the $s$
variables, evaluated at
$s_1=1$, $s_i=0$ for $i\neq 1$, can be expressed as linear
combinations of derivatives with respect to the $t$ variables,
evaluated at the same point. This had been previously conjectured,
and proved in a few special cases, by Witten [17]. Our idea is
that it is precisely this result which, when interpreted
geometrically, should provide the sought-for link between
combinatorial and algebro-geometric classes. In fact, this should be
the case even on the Deligne-Mumford compactification
$\mgnbar{g}{n}$. Why we believe this is explained in detail in
section 3. We can show two things. First of all that our idea works
in complex codimension 1 (and we are pretty sure it also works in
codimension two). This is the content of section 4. Secondly that in
all the cases when we have been able to make the Di Francesco,
Itzykson and Zuber correpondence explicit (the first 11 cases,
according to weight, as defined in section 3), this correspondence
translates into identities of the type
$$
\int_{\wmn{m_*}{n}}\prod\psi_i^{d_i}
=\int_{\mgnbar{g}{n}}X_{m_*,n}\prod\psi_i^{d_i}\,,
$$
where the $X_{m_*,n}$ are explicit polynomials in the
algebro-geometric classes, for any choice of $d_1,\dots,d_n$. In
other words, as linear functionals on a large subspace of the
cohomology group of complementary degree, the classes
$X_{m_*,n}$ behave as duals of the cycles $\wmn{m_*}{n}$. The
codimension one case is settled precisely by showing that this
subspace is as large as it can possibly be, as soon as
$n>1$.
\medskip We are grateful to Edward Witten for useful
conversations, and to Enrico Bombieri and Giulia Galbiati for help
with the computational aspects of the present work. Finally we
would like to thank the Institute for Advanced Study, and in
particular Enrico Bombieri and Phil Griffiths, for a very warm
hospitality.

\vskip1.5 truecm
\specialhead 1. Mumford classes
\endspecialhead
As is customary, we denote by $\mgnbar{g}{n}$
the moduli space of stable $n$-pointed genus $g$ algebraic curves and by
$\mgn{g}{n}$ the moduli space of smooth ones. We shall consistently view
these moduli spaces as orbifolds;
likewise, morphisms between moduli spaces will always be morphism of orbifolds.
We denote by $\Sigma_k$ the subspace of $\mgnbar{g}{n}$ consisting of all
stable
curves with exactly $k$ singular points and their specializations; the
codimension of
$\Sigma_k$ is $k$. This gives a stratification of $\mgnbar{g}{n}$ whose
codimension one stratum $\Sigma_1$ is nothing but the boundary
$\partial\mgn{g}{n}=\mgnbar{g}{n}\setminus \mgn{g}{n}$. To describe the
components of the $\Sigma_k$ it is convenient to proceed as follows. Let $C$
be a stable $n$-pointed genus $g$ curve. Denote by $\Gamma$ the
dual graph of $C$. This is the directed graph whose vertices are the
irreducible components of $C$ and whose edges are the nodes of $C$,
equipped with two additional data.
If we let $V=V_\Gamma$ be the set of its vertices,
the first is the map $p=p_\Gamma$
from $V$ to the non-negative integers assigning to each vertex $v$ the genus
$p_v=p_{\Gamma,v}$ of the normalization of the corresponding component of $C$.
The second is a partition of $\{1,\dots,n\}$ indexed by $V$, that is, a map
$P=P_{\Gamma}:V\to\Cal P(\{1,\dots,n\})$ such that $\{1,\dots,n\}$ is the
disjoint
union of the $P(v)$, for $v\in V$. The partition in question is defined as
follows:
$P(v)$ is the set of all indices $i\in \{1,\dots,n\}$ such that the $i$-th
marked
point belongs to the component corresponding to $v$. We set $h_v=\# P(v)$, and
denote by $l_v$ the valency of $v$, that is, the cardinality of the set $L_v$
of
half-edges issuing from $v$. Clearly, the following hold
$$
\eqalign
{g&=\sum_{v\in V}p_v +h^1(\Gamma)\,,\cr
n&=\sum_{v\in V}h_v\,.
}\leqno{(1.1)}
$$
Also, the stability of $C$ translates into
$$
2p_v - 2 + h_v + l_v > 0\leqno{(1.2)}
$$
for every $v\in V$.
It is clear that any $\Gamma$ satisfying (1.1)
and (1.2) arises as the dual graph of a stable $n$-pointed genus $g$ curve. Now
fix a $\Gamma$ as above, and in addition
choose, for each $v\in V$, an ordering of $L_v$. This determines a
morphism
$$
\xi_{\Gamma}: \prod_{v\in
V}\mgnbar{p_v}{h_v+l_v}\to\mgnbar{g}{n}\,,\leqno{(1.3)}
$$
which is defined as follows. A point in the domain is the assignment of an
$(h_v+l_v)$-pointed curve $(C_v; x_{v,1},\dots,x_{v,h_v+l_v})$ of genus $p_v$
for each $v\in V$. The image point under $\xi_{\Gamma}$ is the stable
$n$-pointed curve of genus $g$ that one obtains by identifying two marked
points $x_{v,h_v+i}$ and $x_{w,h_w+j}$ whenever the $i$-th element of $L_v$
and the $j$-th element of $L_w$ are the two halves of an edge; the marked
points of this curve are the images of the points $x_{v,i}$ for $v\in V$ and
$1\leq i\leq h_v$, with the ordering induced by $P$.

The morphism $\xi_{\Gamma}$ is a finite map onto an irreducible component of
$\Sigma_k$, where $k$ is the number of edges of $\Gamma$. As the reader
may easily verify, this component does not depend on the choice of orderings
of the $L_v$. This is a partial justification for omitting mention of these
orderings in the notation for the map $\xi_{\Gamma}$. More importantly,
our reason for introducing the morphisms $\xi$ is to be able to describe
(boundary) cohomology classes on $\mgnbar{g}{n}$ as pushforwards of classes
on $\prod_{v\in V}\mgnbar{p_v}{h_v+l_v}$, and it will turn out that the classes
we
shall so obtain will always be independent of the choice of orderings of the
$L_v$.

Let us denote by $\Delta_{\Gamma}$ the
image of $\xi_{\Gamma}$. The degree of $\xi_{\Gamma}$, as a map to
$\Delta_{\Gamma}$, is precisely $\#\operatorname{Aut}(\Gamma)$.
This has to be taken with a grain of salt, i. e., is true only if one regards
$\xi_{\Gamma}$ as a morphism of orbifolds. For example, given the graph
$\Gamma$ in Figure 1a), the corresponding map
$$
\xi: \mgnbar{0}{3}\times\mgnbar{0}{3}\to\Delta_{\Gamma}\subset
\mgnbar{2}{0}
$$
is set-theoretically one-to-one, since both source and target consist of one
point; on the other hand, while the unique point of
$\mgnbar{0}{3}\times\mgnbar{0}{3}$ has no automorphisms, the
automorphism group of its image in $\Delta_{\Gamma}$ has order twelve,
as the automorphism group of $\Gamma$; thus the degree of $\xi$ equals twelve
in this case.

It is clear that the
components of $\Sigma_k$ are precisely the $\Delta_{\Gamma}$ for
which $\Gamma$ has $k$ edges. Moreover $\Delta_{\Gamma}$
and $\Delta_{\Gamma'}$ are equal if and only if $\Gamma$
and $\Gamma'$ are isomorphic. We shall write $\delta_\Gamma$ to denote the
fundamental class of $\Delta_\Gamma$ in the rational cohomology of
$\mgnbar{g}{n}$.

To exemplify, let's look at the codimension 1 case. The possible dual graphs
are
illustrated (with repetitions) in Figure 1b). The meaning of the
labelling should be fairly clear; for instance, in the graph $\Gamma_{q,I}$ the
function $p$ assigns the integers $q$ and $g-q$ to the two vertices, as
indicated,
while the partition $P_{\Gamma_{q,I}}$ is just $\{I,\Cal CI\}$.

To abbreviate, we
shall write $\xi_{irr}$, $\xi_{q,I}$, $\delta_{irr}$, $\delta_{q,I}$, instead
of
$\xi_{\Gamma_{irr}}$, $\xi_{\Gamma_{q,I}}$, ad so on. In addition, we set
$\delta=\sum \delta_{\Gamma}$, where $\Gamma$ runs through all {\it
isomorphism classes} of codimension 1 dual graphs. For each of the dual graphs
illustrated in Figure 1b), the underlying directed graph has an order two
automorphism. This induces an automorphism of $\Gamma_{irr}$ and an
isomorphism between $\Gamma_{q,I}$ and $\Gamma_{g-q,\Cal CI}$ (which,
incidentally, is an automorphism precisely when $q=g/2$ and $n=0$). It follows
that
$$
\delta={1\over 2}{\xi_{irr}}_*(1) + {1\over 2}\hskip-0.3truecm
\sum_{0\leq q\leq g\atop
I\subset
\{1,\dots,n\}}\hskip-0.3truecm{\xi_{q,I}}_*(1)\,.\leqno{(1.4)}
$$
Notice that the pushforwards in this formula are well defined since Poincar\'e
duality with rational coefficients holds for all the moduli spaces involved.
Next we
look at the morphism ``forgetting the last marked point''
$$
\pi_{n+1}:\mgnbar{g}{n+1}\longrightarrow \mgnbar{g}{n}\,,
$$
which we may also view as the universal curve $\overline{\Cal C}_{g,n}\to
\mgnbar{g}{n}$. We denote by $\sigma_1,\dots,\sigma_n$ the canonical sections
of
$\pi_{n+1}$, and by $D_1,\dots,D_n$ the divisors in $\mgnbar{g}{n+1}$ they
define.
A point of $D_i$ corresponds to a stable $(n+1)$-pointed genus $g$ curve
obtained
by attaching at one point a $\Bbb P^1$ to a curve $C$ of genus $g$, with the
$i$-th and $(n+1)$-th marked points on the $\Bbb P^1$, and the remaining ones
on $C$. More exactly, $D_i$ is just $\Delta_{\Gamma}$, where $\Gamma$ is
the graph with one edge and two vertices $\{v,w\}$, $p_v=0$, $p_w=g$,
$P(v)=\{i,n+1\}$, and $P(w)=\{1,\dots,i-1,\widehat{i},i+1,\dots,n\}$.
We let $\omega_{\pi_{n+1}}$ be the relative dualizing sheaf. We set
$$
\eqalign{\psi_i&=c_1(\sigma_i^*(\omega_{\pi_{n+1}}))\,,\cr
K&=c_1\left(\omega_{\pi_{n+1}}\left(\sum D_i\right)\right)\,,\cr
\kappa_i&={\pi_{n+1}}_*(K^{i+1})\,.
}
$$
Here, of course, Chern classes are taken to be in rational cohomology, and
${\pi_{n+1}}_*$ is well defined since, as we already observed, Poincar\'e
duality
holds for both the domain and the target of $\pi_{n+1}$. We shall call the
classes
$\kappa_i$ {\it Mumford classes}; in fact, for
$n=0$, their analogues in the intersection ring were first introduced by
Mumford in [12]. Notice that
$\kappa_0=2g-2+n$. Of course, a possible alternative generalization of
Mumford's
$\kappa$'s to the case of
$n$-pointed curves would be the classes
$$
\tilde{\kappa}_i={\pi_{n+1}}_*(c_1(\omega_{\pi_{n+1}})^{i+1})\,.
$$
However, these are not as nicely behaved, from a functorial point of view, as
the $\kappa$'s, as we shall presently see. At any rate, the two are related by
$$
\kappa_a=\tilde{\kappa}_a+\sum_{i=1}^n\psi_i^a\,.\leqno{(1.5)}
$$
The proof of this formula is based on the observation that, for any $j$, taking
residues along
$D_j$ gives an isomorphism between the restriction of
$\omega_{\pi_{n+1}}(\sum D_i)$ to $D_j$ and $\Cal O_{D_j}$. For brevity we set
$\pi=\pi_{n+1}$ and $\tilde{K}=c_1(\omega_\pi)$. Then
$$
\eqalign{
\kappa_a&=\pi_*((\tilde{K}+\sum D_i)^{a+1})=\pi_*(\tilde{K}^{a+1})
+\sum_{l=0}^a\sum_{i=1}^n{a+1\choose l}\pi_*(\tilde{K}^l\cdot D_i^{a-l+1})\cr
&=\pi_*(\tilde{K}^{a+1}) +\sum_{l=0}^a\sum_{i=1}^n{a+1\choose
l}\pi_*\left({\tilde{K}^l\big|_{D_i}}\cdot
(-\tilde{K})^{a-l}\big|_{D_i}\right)\cr
&=\pi_*(\tilde{K}^{a+1}) +\sum_{l=0}^a(-1)^{a-l}{a+1\choose
l}\cdot\sum_{i=1}^n\psi_i^a\cr
&=\tilde{\kappa}_a +\sum_{i=1}^n\psi_i^a\,.
}
$$
Some of the good properties enjoyed by the $\kappa_i$ but not by the
$\tilde{\kappa}_i$ are
$$
\kappa_1\ \text{ is ample}\,,\leqno{(1.6)}
$$
$$
{\pi_n}_*(\psi_1^{a_1}\dots\psi_{n-1}^{a_{n-1}}\psi_n^{a_n+1})=
\psi_1^{a_1}\dots\psi_{n-1}^{a_{n-1}}\kappa_{a_n}\,, \leqno{(1.7)}
$$
$$
\xi_{\Gamma}^*(\kappa_a)=\sum_{v\in V}pr_v^*(\kappa_a)\,,\leqno{(1.8)}
$$
where $\xi_{\Gamma}$ is as in (1.3) and $pr_v$ stands for the projection
from $\prod_{v\in V}\mgnbar{p_v}{h_v+l_v}$ onto its $v$-th factor. That (1.6)
holds is fairly well known; a short proof can be found in [2]. In view of
(1.5), formula (1.7) can also be written in the form
$$ {\pi_n}_*(\psi_1^{a_1}\dots\psi_{n-1}^{a_{n-1}}\psi_n^{a_n+1})=
\psi_1^{a_1}\dots\psi_{n-1}^{a_{n-1}}\tilde{\kappa}_{a_n}+
\sum_{j=1}^{n-1}\psi_1^{a_1}\dots\psi_j^{a_j+a_n}\dots\psi_{n-1}^{a_{n-1}}\,.
$$ As such, with the obvious changes in notation, it is part of formula (1) in
[8]; incidentally, the other part is the so-called string equation
$$ {\pi_n}_*(\psi_1^{a_1}\dots\psi_{n-1}^{a_{n-1}})=
\sum_{\{j:a_j>0\}}\psi_1^{a_1}\dots\psi_j^{a_j-1}\dots\psi_{n-1}^{a_{n-1}}\,.
\leqno{(1.9)}
$$ The proof of formula (1) of [8] is essentially given by Witten in section
2b) of [16]. We now come to (1.8). Consider the diagram
$$
\CD
\Cal Y  @>\eta>>  \overline{\Cal C}_{g,n}\\ @V{\pi'}VV  @V{\pi_{n+1}}VV\\
\dsize{\prod_{v\in V}}\mgnbar{p_v}{h_v+l_v}  @>{\xi_\Gamma}>> \mgnbar{g}{n}
\endCD
$$
where
$$
\Cal Y=\coprod_{v\in V}\left(\overline{\Cal C}_{p_v,h_v+l_v}\times
\prod_{u\in V\atop u\neq v}\mgnbar{p_u}{h_u+l_u}\right)
$$
and $\eta$ is defined by glueing along sections in the manner prescribed by
$\Gamma$. We also let $\pi_v'$ and $\eta_v$ be the restrictions of $\pi'$ and
$\eta$ to
$\overline{\Cal C}_{p_v,h_v+l_v}\times \prod_{u\neq v}\mgnbar{p_u}{h_u+l_u}$.
The morphism $\pi_v'$ is endowed with $h_v+l_v$ canonical sections
$S_1,\dots,S_{h_v+l_v}$. Now the point is that, by the very definition of
dualizing
sheaf, $\eta_v^*(K)=c_1(\omega_{\pi_v'}(\sum S_i))$. Property (1.8) follows.

Another useful property of the classes $\kappa$ is that, on $\mgnbar{g}{n}$,
one
has
$$
\kappa_a=\pi_n^*(\kappa_a)+\psi_n^a\,.\leqno{(1.10)}
$$
To prove this, look at the diagram
$$
\CD
\overline{\Cal C}_{g,n}  @>{\lambda}>>  \overline{\Cal C}_{g,n-1}\\
@V{\varphi_{n+1}}VV  @V{\varphi_{n}}VV\\
\mgnbar{g}{n}  @>{\pi_n}>> \mgnbar{g}{n-1}
\endCD
$$
and denote by $D_1,\dots,D_n$ (resp.,
$D_1',\dots,D_{n-1}'$) the canonical sections of $\varphi_{n+1}$ (resp.,
$\varphi_{n}$). We claim that
$$
\lambda^*\left(\omega_{\varphi_n}\left(\sum
D_i'\right)\right)=\omega_{\varphi_{n+1}}\left(\sum_{i<n}
D_i\right)\,.\leqno{(1.11)}
$$
In fact, there is a natural homomorphism from
$\lambda^*(\omega_{\varphi_n}(\sum D_i'))$ to $\omega_{\varphi_{n+1}}(\sum
D_i)$; we wish to see that this is an isomorphism onto
$\omega_{\varphi_{n+1}}(\sum_{i<n} D_i)$. The question is local in the orbifold
sense. It is therefore sufficient to prove (1.11) when universal curves over
moduli
are replaced by Kuranishi families. To keep things simple we shall use the same
notation in this new setup. Therefore, from now on, $\varphi_n=\pi_n:\Cal C\to
B$
will stand for a Kuranishi family of stable $(n-1)$-pointed curves, and
$D'_1,\dots,D'_{n-1}$ for its canonical sections. A suitable blow-up
$\Cal C'$ of $\Cal C\times_B\Cal C$ provides a Kuranishi family
$\varphi_{n+1}:\Cal C'\to \Cal C$ of stable $n$-pointed curves, whose canonical
sections we shall denote by
$D_1,\dots,D_n$. The diagram we shall look at is
$$
\CD
\Cal C'  @>{\lambda}>>  \Cal C\\ @V{\varphi_{n+1}}VV  @V{\varphi_{n}}VV\\
\Cal C  @>{\pi_n}>> B
\endCD
$$
Observe that, in order to prove that the natural homomorphism from
$\lambda^*(\omega_{\varphi_n}(\sum D_i'))$ to $\omega_{\varphi_{n+1}}(\sum
D_i)$ is an isomorphism onto $\omega_{\varphi_{n+1}}(\sum_{i<n} D_i)$,
it suffices to do so fiber by fiber.
Now  a fiber of $\varphi_n$ is nothing but an $(n-1)$-pointed curve
$(C;x_1,\dots,x_{n-1})$. The inverse image via $\lambda$ of this fiber, which
we
denote by $X$, can be described as follows. Set $Y=C\times C$, $D_i''=C\times
\{x_i\}$ for $i<n$, and let $D_n''$ be the diagonal. Then $X$ is the blow-up of
$Y$ at
the points where $D_n''$ hits a node or one of the $D_i''$, for $i<n$. We
denote the
exceptional curves arising from points of this second type by
$E_1,\dots,E_{n-1}$.
We also observe that the proper transform of each $D_i''$ is just the
intersection of
$D_i$ with $X$. But then
$$
\lambda^*(\omega_{\varphi_n}({\textstyle\sum} D_i'))\big|_X=
\lambda^*(\omega_{\varphi_n})({\textstyle\sum}_{i<n}
D_i)\big|_X({\textstyle\sum} E_i)\,.
$$
On the other hand it is easy to prove (cf. [16]) that
$\lambda^*(\omega_{\varphi_n})=
\omega_{\varphi_{n+1}}(-\Delta)$, where $\Delta$ is the divisor in $\Cal
C'$ defined as follows. Look at fibers of $\varphi_{n+1}$ containing a smooth
rational component meeting the rest of the fiber at only one point; we will
refer to
such a component as a ``rational tail''. Then
$\Delta$ is the divisor swept out by the rational tails containing only two
marked
points, one of which is the $n$-th point. Notice that the divisor cut out by
$\Delta$
on
$X$ is $\sum E_i$. Coupled with the formula above, this implies that
$$
\lambda^*(\omega_{\varphi_n}({\textstyle\sum} D_i'))\big|_X=
\omega_{\varphi_{n+1}}({\textstyle\sum}_{i<n} D_i)\big|_X\,,
$$
which is what we had to show. To prove (1.10) we now argue exactly as in the
proof of (1.5), but applying ${\varphi_{n+1}}_*$ to the $(a+1)$-st
selfintersection of
both sides of (1.11) instead of using the definition of $\kappa_a$.

Using (1.10) and the push-pull formula, we can apply formula (1.7) repeatedly
to
obtain
$$
\eqalign{
(\pi_{n-1}\pi_n)_*(\psi_1^{a_1}\dots\psi_{n-2}^{a_{n-2}}\psi_{n-1}^{a_{n-1}+1}
\psi_n^{a_n+1})&={\pi_{n-1}}_*(\psi_1^{a_1}\dots\psi_{n-2}^{a_{n-2}}
\psi_{n-1}^{a_{n-1}+1}\kappa_{a_n})\cr
&={\pi_{n-1}}_*(\psi_1^{a_1}\dots\psi_{n-1}^{a_{n-1}+1}({\pi_{n-1}}^*
(\kappa_{a_n})+
\psi_{n-1}^{a_{n}}))\cr
&=\psi_1^{a_1}\dots\psi_{n-2}^{a_{n-2}}(\kappa_{a_{n-1}}\kappa_{a_n}+
\kappa_{a_{n-1}+a_n})\,, }
$$
$$
\eqalign{ (\pi_{n-2}\pi_{n-1}\pi_n)_* (\psi_1^{a_1}\dots\psi_{n-3}^{a_{n-3}}&
\psi_{n-2}^{a_{n-2}+1}\psi_{n-1}^{a_{n-1}+1}\psi_n^{a_n+1})=
(\psi_1^{a_1}\dots\psi_{n-3}^{a_{n-3}})\times\hskip3truecm\cr
\times(\kappa_{a_{n-2}}\kappa_{a_{n-1}}\kappa_{a_n}+
\kappa_{a_{n-2}}\kappa_{a_{n-1}+a_n}&+
\kappa_{a_{n-1}}\kappa_{a_{n-2}+a_n} +\kappa_{a_n}\kappa_{a_{n-2}+a_{n-1}}+
2\kappa_{a_{n-2}+a_{n-1}+a_n})\,,
}
$$
and so on. In general, one finds formulas
$$
(\pi_{k+1}\dots\pi_n)_*(\psi_1^{a_1}\dots\psi_k^{a_k}
\psi_{k+1}^{a_{k+1}+1}\dots\psi_n^{a_n+1})=\psi_1^{a_1}\dots\psi_k^{a_k}
R_{a_{k+1}\dots a_n}\,,
\leqno{(1.12)}
$$
where $R_{a_{k+1}\dots a_n}$ is a polynomial in the Mumford classes. A
compact expression for $R_{b_1\dots b_l}$, which we learned from Carel Faber,
is
$$
R_{b_1\dots b_l}=\sum_{\sigma\in\Cal
S_l}\hskip-0.2truecm\kappa_\sigma\,,
\leqno{(1.13)}
$$
where $\kappa_\sigma$ is defined as follows. Write the permutation $\sigma$
as a product of $\nu(\sigma)$ disjoint cycles, including 1-cycles:
$\sigma=\alpha_1\dots\alpha_{\nu(\sigma)}$, where we think of $\Cal S_l$ as
acting on the $l$-tuple $(b_1,\dots ,b_l)$. Denote by $|\alpha|$ the sum of
the elements of a cycle $\alpha$. Then
$$
\kappa_\sigma=\kappa_{|\alpha_1|}\kappa_{|\alpha_2|}\dots
\kappa_{|\alpha_{\nu(\sigma)}|}\,.
$$
Formula (1.12), together with the string equation, expresses the remarkable
fact
that the intersection theory of the classes $\psi_i$ and $\kappa_i$ on a fixed
$\mgnbar{g}{n}$ is completely determined by the intersection theory of the
$\psi_i$ alone on all the $\mgnbar{g}{\nu}$ with $\nu\geq n$, and conversely.
A special case of this is Witten's remark [16] that knowing the intersection
numbers of the $\kappa$'s on $\mgnbar{g}{0}$ is equivalent to knowing the
intersection numbers of the $\psi$'s on all the $\mgnbar{g}{n}$. Using the
``correct'' classes $\kappa_i$ makes all of this particularly transparent.
\vskip1truecm
A final remark has to do with Wolpert's formula [18] stating that,
on $\mgnbar{g}{0}$,
$$
\kappa_1={1\over 2\pi^2}[WP]\,,
$$
where $WP$ is the Weil-Petersson K\"ahler form. It may be observed that this
carries over with no formal changes to $\mgnbar{g}{n}$, for any $n$. To
prove the formula (including the case considered by Wolpert), one may proceed
as
follows. The ``restriction phenomenon'' (page 502 of Wolpert's paper) amounts
to
the statement that the analogue of (1.8) above holds for the class of the
Weil-Petersson K\"ahler form. Arguing by induction on the genus and the
number of marked points, we may then assume that the difference between
$\kappa_1$ and ${1\over 2\pi^2}[WP]$ restricts to zero on any component of
the boundary of $\mgnbar{g}{n}$. One then proves a general lemma to the
effect that a degree two cohomology class with this property actually vanishes
on $\mgnbar{g}{n}$, except in the cases when $\mgnbar{g}{n}$ is
one-dimensional; these are the initial cases of the induction and are dealt
with
by direct computation. The general lemma is proved, although not formally
stated, in [1], for $n=0$; similar ideas can be used to deal with the case when
$n>0$.

\vskip1.5truecm
\specialhead 2. Combinatorial classes
\endspecialhead
Following Kontsevich [7], whose notation we shall adhere to throughout this
section, we consider connected ribbon graphs with metric and with valency of
each vertex greater than or equal to three such that the corresponding
noncompact surface has genus $g$ and $n$ punctures, numbered by
$\{1,\dots,n\}$. We let $\mcomb{g}{n}$ be the space of equivalence classes of
such graphs, endowed with its natural orbifold structure. Recall that the map
$$
\mgn{g}{n}\times{\Bbb R}^n_+ \to \mcomb{g}{n}
$$
which associates to a smooth $n$-punctured curve and an $n$-tuple of
positive real numbers the critical graph of the corresponding canonical Strebel
quadratic differential is a homeomorphism or orbifolds. As Kontsevich has
indicated, the above map extends to a map of ``orbispaces''
$$
\mgnbar{g}{n}\times{\Bbb R}^n_+ \to \mcombar{g}{n}\,,
$$
where $\mcombar{g}{n}$ is a suitable partial compactification of
$\mcomb{g}{n}$. This map, however, is no longer one-to-one, as a certain
amount of contraction takes place at the boundary. More specifically,
$\mcombar{g}{n}$ is isomorphic to $\mprime{g}{n}\times{\Bbb R}^n_+$, where
$\mprime{g}{n}$ equals $\mgnbar{g}{n}$ modulo the {\it closure} of the
following equivalence relation. Two stable $n$-pointed curves are considered
equivalent if there is a homeomorphism of pointed curves between them which
is complex analytic on all components containing at least one marked point.
We let
$$
\alpha: \mgnbar{g}{n}\to\mprime{g}{n}
$$
be the natural projection.

Now fix a sequence $m_*=(m_0,m_1,\dots)$ of non-negative integers almost
all of which are zero. We denote by $\mgn{m_*}{n}$ the space of equivalence
classes of connected numbered ribbon graphs with metric having $n$ boundary
components, $m_i$ vertices of valency $2i+1$ for each $i$, and no vertices of
even
valency. The dimension of $\mgn{m_*}{n}$ is nothing but the number of edges of
such a
graph, and hence
$$
\dim_{\Bbb R}\mgn{m_*}{n}={1\over 2}\sum_i m_i(2i+1)\,.
$$
When $m_0=0$, the space $\mgn{m_*}{n}$ naturally lies inside $\mcomb{g}{n}$,
where $g$ is given by the formula
$$
2g-2+n={1\over 2}\sum_i m_i(2i-1)\,.
$$
More generally, the Strebel construction always gives a map from
$\mgn{m_*}{n}$ to $\mgn{g}{n}\times \Bbb R_+^n$, even for $m_0\not=0$, so
that in
particular the classes $\psi_i$ can be pulled back to $\mgn{m_*}{n}$. In all
cases
we have
$$
\dim_{\Bbb R}\mgn{m_*}{n}=6g-6+2n-\sum_i m_i(i-1)=\dim_{\Bbb R}
\mgn{g}{n}- \sum_im_i(i-1)\,.
$$
On each component of $\mgn{m_*}{n}$ one can put a natural orientation, as
explained
on page 11 of [7].
When $m_0=0$ it can be seen that, with this orientation, $\mgn{m_*}{n}$ is a
cycle with
non-compact support in
$\mcombar{g}{n}$. As such, it defines a class
$$
[\mgn{m_*}{n}]\in H^{non-compact}_{d+n-2k}(\mcombar{g}{n},\Bbb Q)\,,
$$
where $d=6g-6+2n$ and $k=\sum_im_i(i-1)$, hence an element of the dual
of
$$
H_c^{d+n-2k}(\mcombar{g}{n},\Bbb Q)=H^{d-2k}(\mprime{g}{n},\Bbb Q)\,.
$$
This can be also viewed as an element
$\wmn{m_*}{n}\in H_{d-2k}(\mprime{g}{n},\Bbb Q)$

It has been conjectured by Kontsevich [7] (and previously, in a somewhat more
restricted form, by Witten) that the classes $\wmn{m_*}{n}$ ``can be
expressed in terms of the Mumford-Miller classes''. We next give a possible
interpretation of this sentence and a more precise form of the conjecture. The
statement
is made a bit clumsy by the fact that it is not a priori clear whether the
classes
$\wmn{m_*}{n}$ lift to classes in $H_*(\mgnbar{g}{n},\Bbb Q)$, as would happen,
for
instance, if Poincar\'e duality held on $\mprime{g}{n}$. What is certainly true
is that,
given a cohomology class
$x\in H^*(\mgnbar{g}{n},\Bbb Q)$,
$$
\varphi\mapsto\int_{\mgnbar{g}{n}}x\cup\alpha^*(\varphi)
$$
defines a linear functional on $H^{d-2k}(\mprime{g}{n},\Bbb Q)$. What may be
conjectured is that this functional equals $\wmn{m_*}{n}$ for an $x$ of the
form
$$
x=P_{m_*,n}(\kappa_1,\kappa_2,\dots)+\beta_{m_*,n}\,,
$$
where $P_{m_*,n}$ is a weighted-homogeneous polynomial in the Mumford
classes and $\beta_{m_*,n}$ is supported on the boundary of moduli. In what
follows we shall often take the liberty of writing $\wmn{m_*}{n}\equiv x$
to express this, when no confusion seems likely.
One may be more precise about $P_{m_*,n}$ and $\beta_{m_*,n}$. Define the {\it
level} of a monomial $\prod_{a\geq 1}\kappa_a^{h_a}$ in the
Mumford classes to be
$\sum_a h_a$.
Then $P_{m_*,n}$ should be of the form
$$
\prod_{i=2}^\infty{(2^i(2i-1)!!)^{m_i}\over m_i!}\prod_{i\geq
2}\kappa_{i-1}^{m_i}+\text{ a linear combination of monomials of lower
level}\,.\leqno{(2.1)}
$$
As for $\beta_{m_*,n}$ it should be a linear combination of classes of the
form ${\xi_{\Gamma}}_*(y)$, where $y$ is a monomial in the Mumford
classes and in the $pr_v^*(\psi_i)$, for $i>h_v$, where of course we have
freely
used the notation established in section 1. As a special case, one should have
$$
\wmn{(0,m_1,0,\dots,0,m_j=1,0,\dots)}{n}\equiv 2^j(2j-1)!!\kappa_{j-1}+\text{
boundary terms}\,.\leqno{(2.2)}
$$
In the next section we shall give evidence for these conjectures. In
section 4 we shall prove (2.2) in the codimension 1 case for $n>1$; more
exactly,
we shall show that
$$
\wmn{(0,m_1,1,\dots)}{n}\equiv 12\kappa_1-\delta\,,\leqno{(2.3)}
$$
where $\delta$ is the usual class of the boundary. We shall also see that, for
$n>1$, this formula includes as a special case the main result of Penner [13],
with the following caveat. In our notation, what Penner claims is that
$\wmn{(0,m_1,1,\dots)}{n}= 6\tilde{\kappa}_1$ on the open moduli space
$\mgn{g}{n}$, while the correct formula is $\wmn{(0,m_1,1,\dots)}{n}=
12\kappa_1$. It should be said that Penner's argument, which, by the way, is
entirely different from ours, is completely correct, except for two minor
mistakes in the interpretation of what has been proved. The first mistake is
that, as we have noticed in section 1, the class of the Weil-Petersson K\"ahler
form is $\kappa_1$ and not $\tilde{\kappa}_1$. The second mistake actually
occurs in Theorem A.2 of [14], where the explicit expression of the
Weil-Petersson K\"ahler form should be divided by two. In fact, if one looks at
how this is
gotten, one sees that it is computed as the pull-back of the Weil-Petersson
K\"ahler form on $\mgn{2g+n-1}{0}$ via the doubling map which associates to
a genus $g$ smooth $n$-pointed curve $C$ the curve obtained by attaching at
the punctures two identical copies of $C$. But now, since one is doubling, one
must also divide by two, as the resulting curve carries the extra automorphism
which exchanges the two components. An advantage of our method
over Penner's is perhaps that, in addition to giving a certain amount of
control
over the boundary, it is not special to the codimension one case but provides,
at least in principle, a mechanism for dealing with classes of higher
codimension.

\vskip1.5 truecm
\specialhead 3. Geometrical consequences of a result of Di Francesco, Itzykson,
and Zuber
\endspecialhead
Following Witten [16] and Kontsevich [7], given a sequence of
non-negative integers
$\ul{d}=(d_1,\dots,d_n)$ and an infinite sequence
$m_*=(m_0,m_1,m_2,\dots)$ of non-negative integers, almost all zero, we set
$$
\langle\tau_{\ul{d}}\rangle_{m_*}=\langle\tau_{d_1}\dots\tau_{d_n}
\rangle_{m_*}=\int_{\mgn{m_*}{n}}\prod_{i=1}^n\psi_i^{d_i}\times[\Bbb R_+^n]\,,
$$
where $[\Bbb R_+^n]$ stands for the fundamental class with compact support of
$\Bbb R_+^n$. This integral is zero unless
$
\sum_{i}d_i= {1\over 2}\dim\mgn{m_*}{n}={1\over 4}\sum_i m_i(2i+1)
$. Notice that, when $m_0=0$, one can also write
$$
\langle\tau_{\ul{d}}\rangle_{m_*}=
\int_{\wmn{m_*}{n}}\prod_{i=1}^n\psi_i^{d_i}\,,
$$
and again this is zero unless
$$
\sum_{i}d_i= {1\over 2}\dim\wmn{m_*}{n}=3g-3+n-\sum_i (i-1)m_i\,,
$$
where $2g-2+n=(1/2)\sum_i m_i(2i-1)$. We also set
$$
\langle\tau_{\ul{d}}\rangle_{g,n}=
\int_{\mgnbar{g}{n}}\prod_{i=1}^n\psi_i^{d_i}\,.
$$
It is clear that
$\langle\tau_{\ul{d}}\rangle_{g,n}=\langle\tau_{\ul{d}}\rangle_{m_*}$ for
$m_*=(0,4g-4+2n,0,0,\dots)$. The symbol $\langle\tau_{\ul{d}}\rangle$, with no
subscripts, stands for $\langle\tau_{\ul{d}}\rangle_{g,n}$ when the number $g$
{\it defined} by $3g-3+n=\sum d_i$ is a non-negative integer, and is set to
zero
otherwise. Sometimes the abbreviated notation
$\tau_0^{n_0}\tau_1^{n_1}\tau_2^{n_2}\dots$ is used in place of
$$
\underbrace{\tau_0\dots\tau_0}_{n_0\text{ times}}
\underbrace{\tau_1\dots\tau_1}_{n_1\text{ times}}
\underbrace{\tau_2\dots\tau_2}_{n_2\text{ times}}\dots
$$
One then considers the formal power series
$$
\eqalign{ F(t_*,s_*)&=\sum_{n_*,m_*}\langle\tau_{\ul{d}}\rangle_{m_*}
{t_*^{n_*}\over n_*!}s_*^{m_*}\,,\cr Z(t_*,s_*)&=\exp(F(t_*,s_*))\,,
}
$$
where the
following notational conventions are adopted. First of all
$$
t_*=(t_0,t_1,t_2,\dots)\,,\qquad s_*=(s_0,s_1,s_2,\dots)\,,
$$
are infinite sequences of indeterminates, and
$$
m_*=(m_0,m_1,m_2,\dots)\,,\qquad n_*=(n_0,n_1,n_2,\dots)\,,
$$
are infinite sequences of non-negative integers, almost all zero. We have also
set
$$
n_*!=\prod_{i=0}^\infty n_i!\,,\qquad t_*^{n_*}=\prod_{i=0}^\infty t_i^{n_i}\,,
$$
and similarly for $s_*^{m_*}$. Finally, if $n=\sum_i n_i$, the sequence of
non-negative integers
$\ul{d}=(d_1,\dots,d_n)$ is  determined (up to their order, which is
irrelevant)
by the requirement that $n_i$ equal the number of $j$'s such that $i=d_j$; in
other words, one could also have written
$\langle\prod_{i=0}^\infty\tau_i^{n_i}\rangle_{m_*}$ instead of
$\langle\tau_{\ul{d}}\rangle_{m_*}$.

Now let $\Cal{H}_N$ be the space of $N\times N$ hermitian matrices, and
consider on it the $U(N)$-invariant measure
$$
dX=\prod_{1\leq i\leq N}dX_{ii}\prod_{1\leq i<j\leq
N}d\operatorname{Re}X_{ij}\,d\operatorname{Im}X_{ij}\,.
$$
For any positive definite $N\times N$ diagonal matrix $\Lambda$ we also
consider the measure
$$
d\mu_\Lambda=c_{\Lambda,N}\exp\left(-{1\over
2}\operatorname{tr}(X^2\Lambda)\right)dX\,,
$$
where $c_{\Lambda,N}$ is the constant such that $\int d\mu_\Lambda=1$. It
has been shown by Kontsevich that, for any fixed $s_*$, and with the
substitution
$$
t_i=-(2i-1)!!\trace(\Lambda^{-2i-1})\,,
$$
the series $Z(t_*,s_*)$ is an asymptotic expansion of the integral
$$
\int_{\Cal{H}_N}\exp\left(-\sqrt{-1}\sum_{j=0}^\infty\left(-{1\over
2}\right)^j s_j {\trace(X^{2j+1})\over
2j+1}\right)d\mu_\Lambda \,\leqno{(3.1)}
$$
as $\Lambda^{-1}$ goes to zero (notice the minus sign in front of the
argument of the exponential, which is missing in the formula given in [7]). To
simplify notations, we set
$$
\eqalign{\br{f}&=\int_{\Cal{H}_N} f d\mu_\Lambda\,,\cr
\brr{f}&=\int_{\Cal{H}_N} f\exp\left({\sqrt{-1}\trace X^3\over 6}\right)
d\mu_\Lambda\,.}
$$
Let us now fix non-negative integers $m_2, m_3, \dots$, almost all equal to
zero, and
set $\hat{s}_*=(0,1,0,0,0\dots)$. It follows from the definitions that
$$
\prod_{i\geq 2}{1\over m_i!}\left({\partial\over\partial
s_i}\right)^{m_i}\hskip-0.2truecm F(t_*,s_*)\
\big|_{s_*=\hat{s}_*}=\sum_{n_*,m_1}
\langle\tau_{\ul{d}}\rangle_{(0, m_1, m_2, m_3, \dots)} {t_*^{n_*}\over
n_*!}\,.\leqno{(3.2)}
$$
In other words, the coefficients of the above derivative of $F$ are just the
intersection numbers
$$
\int_{W_{m_*,n}}\psi_1^{d_1}\dots\psi_n^{d_n}\,,\leqno{(3.3)}
$$
where we have written $m_*$ for $(0, m_1, m_2, m_3, \dots)$. If the conjecture
formulated in section 2 holds true, it should be possible to write these
intersection
numbers under the form
$$
\int_{\mgnbar{g}{n}}(P_{m_*,n}(\kappa_1,\kappa_2,\dots)+\beta_{m_*,n})
\psi_1^{d_1}\dots\psi_n^{d_n}\,,\leqno{(3.4)}
$$
where $P_{m_*,n}$ and $\beta_{m_*,n}$ are as in that section. We contend that
this
result should be implicitly contained in a theorem, conjectured by Witten, and
proved
by Di Francesco, Itzykson, and Zuber [3]. To explain this, the first step is to
observe that, by differentiating (3.1), we obtain asymptotic expansions
$$
\bbrr{\prod_i\left(-\sqrt{-1}\left(-1\over 2\right)^i{\trace X^{2i+1}\over
2i+1}\right)^{\nu_i}}\sim\ \prod_i\left({\partial\over\partial s_i}
\right)^{\nu_i}\hskip-0.2truecm Z(t_*,s_*)\ \big|_{s_*=\hat{s}_*}\,,
$$
for any sequence $\nu_*=(\nu_0,\nu_1,\dots)$ of non-negative integers such that
$\nu_i=0$ for large enough $i$. Now the theorem of Di Francesco, Itzykson, and
Zuber
(henceforth referred to as the DFIZ theorem) states that, given any polynomial
$Q$ in
the {\it odd} traces of
$X$, there exists a differential polynomial
$R_Q=R_Q\left({\partial\over\partial
t_0},{\partial\over\partial t_1},\dots\right)$ such that
$$
\brr{Q} = R_Q Z(t_*)\,,
$$
where $Z(t_*)$ stands for $Z(t_*,\hat{s}_*)$. Putting this together with the
previous
remark shows that
$$
\prod_i\left({\partial\over\partial s_i}\right)^{\nu_i}\hskip-0.2truecm
Z(t_*,s_*)\
\big|_{s_*=\hat{s}_*}=U_{\nu_*}\left({\partial\over\partial t_0},
{\partial\over\partial t_1},\dots\right)Z(t_*)\,,
$$
where $U_{\nu_*}$ is a polynomial. In terms of $F$, this amounts to saying that
$$
\prod_i\left({\partial\over\partial s_i}\right)^{\nu_i}\hskip-0.2truecm
F(t_*,s_*)\
\big|_{s_*=\hat{s}_*}=\widetilde{U}_{\nu_*}\,,
$$
where $\widetilde{U}_{\nu_*}$ is a polynomial in the partial derivatives of
$F(t_*,s_*)$ with respect to the $t$ variables, evaluated at $s_*=\hat{s}_*$.
The
expression that Di Francesco, Itzykson, and Zuber give for $U_{\nu_*}$, and
hence
implicitly for $\widetilde{U}_{\nu_*}$, is quite complicated. However, if we
define
the {\it weight} of a partial derivative
$$
\prod_i\left({\partial\over\partial t_i}\right)^{\nu_i}\hskip-0.2truecm
F(t_*,s_*)
$$
to be $\sum(2i+1)\nu_i$, and the weight of a product of partial derivatives to
be
the sum of the weights of its factors, then what can be said is that
$$
\eqalign{\prod_i\left({\partial\over\partial s_i}\right)^{\nu_i}
&\hskip-0.2truecm
F(t_*,s_*)\ \big|_{s_*=\hat{s}_*}=\cr
\prod_i&
\left(2^i(2i-1)!!{\partial\over \partial t_i}\right)^{\nu_i}\hskip-0.2truecm
F(t_*)+\text{ terms of lower weight}\,,}
\leqno{(3.5)}
$$
where $F(t_*)$ is defined to be equal to $F(t_*,\hat{s}_*)$. In addition, only
terms
whose weight is congruent to $\sum(2i+1)\nu_i$ modulo 3 appear in (3.5).

In the case when $\nu_*=(0,0,m_2,m_3,\dots)$ we have already explained how the
left-hand side of (3.5) is linked to the intersection theory of products of
classes
$\psi_i$ with the $W_{m_*,n}$; it remains to explain the geometric significance
of the
right-hand side. Consider the series
$$
\prod\left({\partial\over\partial t_i}\right)^{\mu_i}\hskip-0.2truecm
F(t_*)=\sum_{n_*}a_{\ul{d}}{t_*^{n_*}\over n_*!}\,.
$$
Then it is easy to show that
$$
a_{\ul{d}}=\left\langle\prod\tau_i^{\mu_i}\tau_{\ul{d}}\right\rangle\,.
$$
In a certain sense one can say that differentiating $F(t_*)$ with respect to
the
$t_i$ variable $\mu_i$ times, for $i=0,1,\dots$, corresponds to the insertion
of
$\prod\tau_i^{\mu_i}$ in the coefficients of $F(t_*)$.
Now fix a positive integer $n$, $\ul{d}=(d_1,\dots,d_n)$, and
$m_*=(0,m_1,m_2,\dots)$. Then
$W_{m_*,n}$ is a cycle in $\mprime{g}{n}$, for a well-determined $g$. Setting
$\nu_*=(0,0,m_2,m_3,\dots)$ and equating coefficients in (3.5) one finds that
$\langle\tau_{\ul{d}}\rangle_{m_*}$ is a linear combination, with rational
coefficients,
of terms of the form
$$
\left\langle\prod_i\tau_i^{\lambda_{i,1}}\tau_{\ul{d}_{I_1}}\right\rangle\dots
\left\langle\prod_i\tau_i^{\lambda_{i,k}}\tau_{\ul{d}_{I_k}}\right\rangle\,,
$$
where $\{I_1,\dots,I_k\}$ is a partition of $\{1,\dots,n\}$ and, for any subset
$I$ of $\{1,\dots,n\}$, we set $\langle\tau_{\ul{d}_I}\rangle=
\langle\prod_{i\in I}\tau_{d_i}\rangle$. Moreover, if we set
$\mu_i=\sum_j\lambda_{i,j}$, then $\sum(2i+1)\mu_i$ is not greater than
$\sum_{i\geq 2}(2i+1)m_i$, and  congruent to it modulo 3. For instance, the
term
coming from the highest weight part of the right-hand side of (3.5) is simply
$$
\eqalign{\prod_{i\geq 2}
(2^i(2i-1)!!)^{m_i}&\left\langle\prod\tau_i^{m_i}\tau_{\ul{d}}
\right\rangle_{g,n+\sum m_i}=\cr &\prod_{i\geq 2}
(2^i(2i-1)!!)^{m_i}\int_{\mgnbar{g}{n}}\left(\prod_{a\geq
1}\kappa_a^{m_{a+1}}+\dots\right)\prod\psi_i^{d_i}\,,}
$$
where we have used the formulas for the pushforwards of products of classes
$\psi_i$ given in section 1. The lower weight terms are considerably more
messy. In
particular there is, a priori, no reason why it should be possible to write
each one of
them under the form
$$
\int_{\mgnbar{g}{n}}\alpha\prod\psi_i^{d_i}\,,
$$
where $\alpha$ is a suitable cohomology class. That this indeed happens, at
least in
all the cases we have been able to compute, depends on some remarkable
cancellations, as we shall presently see. At any rate, we have that
$$
\int_{W_{m_*,n}}\psi_1^{d_1}\dots\psi_n^{d_n}=\left[\int_{\mgnbar{g}{n}}
\left(\prod_{i\geq 2}{1\over
m_i!}(2^i(2i-1)!!\kappa_{i-1})^{m_{i}}\right)\psi_1^{d_1}\dots\psi_n^{d_n}
\right]+\cdots\,,
$$
which can be viewed as a first step in writing the intersection number (3.3) in
the
form (3.4).

To illustrate the procedure we just described we shall work out three examples.
The first one deals with the  cycle $W_{(0, m_1, 1, 0,0,0\dots),n}$.
This  is the only codimension one cycle among the $W_{m_*,n}$ and
corresponds to ribbon graphs having at least one five-valent vertex.
In this case (3.2) reads
$$
{\partial\over\partial s_2}
F(t_*,s_*)\ \big|_{s_*=\hat{s}_*}=\sum_{n_*,m_1}
\langle\tau_{\ul{d}}\rangle_{(0,m_1,1,0,0,\dots)} {t_*^{n_*}\over n_*!}\,.
$$
On the other hand the DFIZ theorem tells us,
in this case, that
$$
{\partial\over\partial s_2}Z(t_*,s_*)\ \big|_{s_*=\hat{s}_*}=
\left(12{\partial\over\partial t_2}-{1\over 2}{\partial^2\over\partial t_0^2}
\right)Z(t_*)\leqno{(3.6)}
$$
As $Z=\exp F$, we then get, upon dividing by $Z$,
$$
{\partial\over\partial s_2}F(t_*,s_*)\ \big|_{s_*=\hat{s}_*}=
12{\partial F\over\partial t_2}(t_*)-
{1\over 2}{\partial^2 F\over\partial t_0^2}(t_*)-{1\over 2}
\left({\partial F\over\partial t_0}(t_*)\right)^2\,.
$$
Comparing coefficients, this amounts to
$$
\eqalign{\langle\tau_{\ul{d}}\rangle_{(0,m_1,1,0,0,\dots)} &=
12\langle\tau_2\tau_{\ul{d}}\rangle_{g,n+1}
-{1\over 2}\langle\tau_0\tau_0\tau_{\ul{d}}\rangle_{g-1,n+2}\cr
&\quad-{1\over 2}\sum_{I\subset\{1,\dots,n\}}
\langle\tau_0\tau_{\ul{d}_I}\rangle_{p,h+1}
\langle\tau_0\tau_{\ul{d}_{\Cal{C}I}}\rangle_{g-p,n-h+1}\,,
}$$
where
$$
\tau_{\ul{d}_I}=\prod_{i\in I}\tau_{d_i}\,,
$$
$g$ is given by
$$
\sum_{i=1}^n d_i=3g-3+n-{1\over 2}\operatorname{codim}
(\wmn{(0,m_1,1,0,0,\dots)}{n})=3g-3+n-1\,,
$$
$m_1$ by
$$
m_1=4g-7+2n\,,
$$
and $h$ and $p$ by
$$
h=\# I\,,\qquad 3p-3+h+1=\sum_{i\in I}d_i\,.
$$
Using (1.7), this can be rewritten as
$$
\eqalign{
\int_{W_{(0,m_1,1,0,0,\dots),n}}\psi_1^{d_1}\dots\psi_n^{d_n}&=
12\int_{\mgnbar{g}{n+1}}\psi_1^{d_1}\dots\psi_n^{d_n}\psi_{n+1}^2
-{1\over
2}\int_{\mgnbar{g-1}{n+2}}\xi_{irr}^*(\psi_1^{d_1}\dots\psi_n^{d_n})\cr
&\quad-{1\over 2}\hskip-0.2truecm\sum_{0\leq p\leq g\atop
I\subset\{1,\dots,n\}}
\int_{\mgnbar{p}{h+1}\times\mgnbar{g-p}{n-h+1}}
\xi_{p,I}^*(\psi_1^{d_1}\dots\psi_n^{d_n})\cr
&=12\int_{\mgnbar{g}{n}}\kappa_1\psi_1^{d_1}\dots\psi_n^{d_n} -{1\over
2}\int_{\mgnbar{g}{n}}{\xi_{irr}}_*(1)\psi_1^{d_1}\dots\psi_n^{d_n}\cr
&\quad-{1\over 2}\hskip-0.2truecm\sum_{0\leq p\leq g\atop
I\subset\{1,\dots,n\}}
\int_{\mgnbar{g}{n}}{\xi_{p,I}}_*(1)\psi_1^{d_1}\dots\psi_n^{d_n}\,, }
$$
or, in view of formula (1.4),
$$
\int_{W_{(0,m_1,1,0,0,\dots),n}}\psi_1^{d_1}\dots\psi_n^{d_n}=
\int_{\mgnbar{g}{n}}(12\kappa_1-\delta)\psi_1^{d_1}\dots\psi_n^{d_n}\,.
\leqno{(3.7)}
$$
Now let's turn to the codimension 2 case. Among the $\wmn{m_*}{n}$ there are
two codimension 2 classes, corresponding to $m_*=(0,m_1,0,1,0,\dots)$ and to
$m_*=(0,m_1,2,0,\dots)$. The first one corresponds to ribbon graphs with at
least one 7-valent vertex, the second one to ribbon graphs with at least two
5-valent vertices. To make notations lighter, from now on we shall adopt the
following convention. Whenever identities between derivatives of $Z$ or $F$
will
be given, these will always be meant to hold for $s_*=\hat{s}_*=(0,1,0,\dots)$,
unless otherwise specified. With this notation, the DFIZ theorem gives
$$
\leqalignno{
{\partial Z\over\partial s_3}&=120{\partial Z\over\partial t_3}-
6{\partial^2 Z\over\partial t_0\partial t_1}
+{5\over 4}{\partial Z\over\partial t_0}\,,&(3.8)\cr
&\phantom{.}\cr
{\partial^2 Z\over\partial s_2^2}&=144{\partial^2 Z\over\partial t_2^2}
-840{\partial Z\over\partial t_3}-
12{\partial^3 Z\over\partial t_0^2\partial t_2}
+24{\partial^2 Z\over\partial t_0\partial t_1}
+{1\over 4}{\partial^4 Z\over\partial t_0^4}
-3{\partial Z\over\partial t_0}\,.&(3.9)
}
$$
In terms of derivatives of $F$ these translate into
$$
{\partial F\over\partial s_3}=120{\partial F\over\partial t_3}
-6{\partial^2 F\over\partial t_0\partial t_1}
-6{\partial F\over\partial t_0}{\partial F\over\partial t_1}
+{5\over 4}{\partial F\over\partial t_0}\leqno{(3.10)}
$$
$$
\eqalign{{\partial^2 F\over\partial s_2^2}+
\left({\partial F\over\partial s_2}\right)^2
&=144{\partial^2 F\over\partial t_2^2}
+144\left({\partial F\over\partial t_2}\right)^2
-840{\partial F\over\partial t_3}
-12{\partial^3 F\over\partial t_0^2\partial t_2}\cr
&\quad-12{\partial^2 F\over\partial t_0^2}{\partial F\over\partial t_2}
-24{\partial F\over\partial t_0}{\partial^2 F\over\partial t_0\partial t_2}
-12\left({\partial F\over\partial t_0}\right)^2{\partial F\over\partial t_2}
+24{\partial^2 F\over\partial t_0\partial t_1}\cr
&\quad+24{\partial F\over\partial t_0}{\partial F\over\partial t_1}
+{1\over 4}{\partial^4 F\over\partial t_0^4}
+{\partial^3 F\over\partial t_0^3}{\partial F\over\partial t_0}
+{3\over 2}{\partial^2 F\over\partial t_0^2} \left({\partial F\over\partial
t_0}\right)^2\cr
&\quad+{1\over 4}\left({\partial F\over\partial t_0}\right)^4
+{3\over 4}\left({\partial^2 F\over\partial t_0^2}\right)^2
-3{\partial F\over\partial t_0}
}
$$
Taking into account (3.6), the second of these yields
$$
\eqalign{
{\partial^2 F\over\partial s_2^2}
&=144{\partial^2 F\over\partial t_2^2}
-840{\partial F\over\partial t_3}
-12{\partial^3 F\over\partial t_0^2\partial t_2}
-24{\partial F\over\partial t_0}{\partial^2 F\over\partial t_0\partial t_2}
+24{\partial^2 F\over\partial t_0\partial t_1}\cr
&\quad+24{\partial F\over\partial t_0}{\partial F\over\partial t_1}
+{1\over 4}{\partial^4 F\over\partial t_0^4}
+{\partial^3 F\over\partial t_0^3}{\partial F\over\partial t_0}
+{\partial^2 F\over\partial t_0^2} \left({\partial F\over\partial t_0}\right)^2
+{1\over 2}\left({\partial^2 F\over\partial t_0^2}\right)^2
-3{\partial F\over\partial t_0}
}\leqno{(3.11)}
$$
As the reader may notice, the right-hand sides of (3.8) and (3.9) contain
considerably fewer terms than one would a priori expect, based on the general
statement of the DFIZ theorem in the form given by (3.5). In fact,
$\partial^4 Z/\partial t_0^4$ is missing from (3.8), while $\partial^7
Z/\partial
t_0^7$, $\partial^3 Z/\partial t_0\partial t_1^2$ and $\partial^5 Z/\partial
t_0^4\partial t_1$ are not present in (3.9). In addition to this phenomenon,
further unexpected cancellations occur when passing from derivatives of $Z$ to
derivatives of $F$. For instance, $(\partial F/\partial t_0)^4$ and $(\partial
F/\partial t_2)^2$ do not appear in (3.11). We shall see in a
moment that these remarkable phenomena have geometrical significance.
Indeed, in equating coefficients in the two sides of (3.10) and (3.11), it is
precisely these facts that make it possible to interpret the resulting
identities as
relations between intersection numbers on a specific moduli space
$\mgnbar{g}{n}$, rather than relations involving intersection numbers on
different moduli spaces.

Term by term, (3.10) translates into
$$
\eqalign{
\langle\tau_{\ul{d}}\rangle_{(0,m_1,0,1,0,\dots)}&=
120\langle\tau_3\tau_{\ul{d}}\rangle_{g,n+1}
-6\langle\tau_0\tau_1\tau_{\ul{d}}\rangle_{g-1,n+2}\cr
&\quad-6\sum_{I\subset\{1,\dots,n\}}
\langle\tau_1\tau_{\ul{d}_I}\rangle_{p,h+1}
\langle\tau_0\tau_{\ul{d}_{\Cal{C}I}}\rangle_{g-p,n-h+1}
+{5\over 4}\langle\tau_0\tau_{\ul{d}}\rangle_{g-1,n+1}\cr
&=120\langle\tau_3\tau_{\ul{d}}\rangle_{g,n}
-6\langle\tau_0\tau_1\tau_{\ul{d}}\rangle_{g-1,n+2}\cr
&\quad-6\sum_{I\subset\{1,\dots,n\}}
\langle\tau_1\tau_{\ul{d}_I}\rangle_{p,h+1}
\langle\tau_0\tau_{\ul{d}_{\Cal{C}I}}\rangle_{g-p,n-h+1}
+30\langle\tau_1\rangle\langle\tau_0\tau_{\ul{d}}\rangle_{g-1,n+1}\,,
}$$
where $3g-3+n-2=\sum d_i$, $m_1=4g-9+2n$, $h=\#I$, $3p-3+h+1=\sum_{i\in
I}d_i$, and we have used the fact that $\langle\tau_1\rangle=1/24$. Proceeding
exactly as in the derivation of (3.7), we conclude that
$$
\int_{W_{(0,m_1,0,1,0,\dots),n}}\psi_1^{d_1}\dots\psi_n^{d_n}=
\int_{\mgnbar{g}{n}}(120\kappa_2+\beta)\psi_1^{d_1}\dots\psi_n^{d_n}\,,
\leqno{(3.12)}
$$
where
$$
\eqalign{
\beta&=-6{\xi_{irr}}_*(\psi_{n+1})-6\hskip-0.3truecm\sum_{0\leq
p\leq g\atop I\subset\{1,\dots,n\}}\hskip-0.3truecm
{\xi_{p,I}}_*(\psi_{h+1}\times 1)\cr
&\quad+30{\xi_{1,\emptyset}}_*(\psi_1\times 1)\,.
}\leqno{(3.13)}
$$
The reader should be warned that $\beta$ is not unambiguously defined, or,
more exactly, that (3.12) holds also for a different choice of boundary term
$\beta$. To see this notice that, using (1.7), we can write
$$
\langle\tau_0\tau_1\tau_{\ul{d}}\rangle_{g-1,n+2}=
(2(g-1)-2+n+1)\langle\tau_0\tau_{\ul{d}}\rangle_{g-1,n+1}=
24(2g+n-3)\langle\tau_1\rangle\langle\tau_0\tau_{\ul{d}}\rangle_{g-1,n+1}\,.
$$
Since $\kappa_0=2g+n-3$ on $\mgnbar{g-1}{n+1}$, this means that we could have
chosen $\beta$ to be given by
$$
\eqalign{
\beta&=-144{\xi_{1,\emptyset}}_*(\psi_{1}\times\kappa_0)-
6\hskip-0.3truecm\sum_{0\leq
p\leq g\atop I\subset\{1,\dots,n\}}\hskip-0.3truecm
{\xi_{p,I}}_*(\psi_{h+1}\times 1)\cr
&\quad+30{\xi_{1,\emptyset}}_*(\psi_1\times 1)\,. }\leqno{(3.13)'}
$$
This kind of ambiguity will be present in all the formulas for combinatorial
classes that we shall give; however, it will be confined to some of the
boundary terms.

Let's turn to formula (3.11). This gives
$$
\eqalign{
2\langle\tau_{\ul{d}}\rangle_{(0,m_1,2,0,0,\dots)}&=
144\langle\tau_2^2\tau_{\ul{d}}\rangle
-840\langle\tau_3\tau_{\ul{d}}\rangle
-12\langle\tau_0^2\tau_2\tau_{\ul{d}}\rangle\cr
&\quad-24\hskip-0.2truecm\sum_{I\sqcup J=\{1,\dots,n\}}\hskip-0.4truecm
\langle\tau_0\tau_{\ul{d}_I}\rangle
\langle\tau_0\tau_2\tau_{\ul{d}_J}\rangle
+24\langle\tau_0\tau_1\tau_{\ul{d}}\rangle\cr
&\quad+24\hskip-0.2truecm\sum_{I\sqcup J=\{1,\dots,n\}}\hskip-0.4truecm
\langle\tau_0\tau_{\ul{d}_I}\rangle
\langle\tau_1\tau_{\ul{d}_J}\rangle
+{1\over 4}\langle\tau_0^4\tau_{\ul{d}}\rangle\cr
&\quad+\hskip-0.1truecm\sum_{I\sqcup J=\{1,\dots,n\}}\hskip-0.4truecm
\langle\tau_0^3\tau_{\ul{d}_I}\rangle
\langle\tau_0\tau_{\ul{d}_J}\rangle\cr
&\quad+\hskip-0.1truecm\sum_{I\sqcup J\sqcup K=\{1,\dots,n\}}
\hskip-0.4truecm
\langle\tau_0^2\tau_{\ul{d}_I}\rangle
\langle\tau_0\tau_{\ul{d}_J}\rangle
\langle\tau_0\tau_{\ul{d}_K}\rangle\cr
&\quad+{1\over 2}\sum_{I\sqcup J=\{1,\dots,n\}}
\hskip-0.4truecm
\langle\tau_0^2\tau_{\ul{d}_I}\rangle
\langle\tau_0^2\tau_{\ul{d}_J}\rangle
-3\langle\tau_0\tau_{\ul{d}}\rangle\,.
}\leqno{(3.14)}
$$
We wish to see that this can be interpreted as an identity among intersection
numbers on $\mgnbar{g}{n}$, where $g$, $m_1$ and the $d_i$ are related by
$\sum d_i=3g-5+n$ and $m_1=4g-10+2n$.
The left-hand side of (3.14) is twice the integral of
$\psi_1^{d_1}\dots\psi_n^{d_n}$
 over $W_{(0,m_1,2,0,0,\dots),n}$. As for the
right-hand side, it is convenient to examine each summand separately. The first
two terms cause no trouble for, using (1.12), they can be written as
$$
\eqalign{144\langle\tau_2\tau_2\tau_{\ul{d}}\rangle_{g,n+2}
-840\langle\tau_3\tau_{\ul{d}}\rangle_{g,n+1}&=
144\int_{\mgnbar{g}{n}}(\kappa_1^2+\kappa_2)
\psi_1^{d_1}\dots\psi_n^{d_n}\cr
&\quad-840\int_{\mgnbar{g}{n}}\kappa_2
\psi_1^{d_1}\dots\psi_n^{d_n}\cr
&=\int_{\mgnbar{g}{n}}(144\kappa_1^2-696\kappa_2)
\psi_1^{d_1}\dots\psi_n^{d_n}\,.
}$$
Now look at the remaining terms.
\medskip\noindent{- Term 3.}
$$
\langle\tau_0^2\tau_2\tau_{\ul{d}}\rangle=
\langle\tau_0^2\tau_2\tau_{\ul{d}}\rangle_{g-1,n+3}=
\int_{\mgnbar{g-1}{n+2}}\kappa_1\prod_{i=1}^n\psi_i^{d_i}=
\int_{\mgnbar{g}{n}}{\xi_{irr}}_*(\kappa_1)\prod\psi_i^{d_i}\,.
$$

\medskip\noindent{- Term 4.} Disregarding the coefficient, this is
$$
\sum_{I\sqcup J=\{1,\dots,n\}}\hskip-0.4truecm
\langle\tau_0\tau_{\ul{d}_I}\rangle_{q,h+1}
\langle\tau_0\tau_2\tau_{\ul{d}_J}\rangle_{r,k+2}\,,
$$
where
$$
h=\#I\,,\quad k=\#J\,,\quad \sum_{i\in I}d_i=3q-3+h+1
\,,\quad \sum_{j\in J}d_j=3r-3+k+2\,.
$$
Since $h+k=n$ and, as we observed above, $\sum d_i=3g-5+n$, this gives
$q+r=g$. Hence
$$
\eqalign{\sum_{I\sqcup J=\{1,\dots,n\}}\hskip-0.4truecm
\langle\tau_0\tau_{\ul{d}_I}\rangle
\langle\tau_0\tau_2\tau_{\ul{d}_J}\rangle
&=\sum_{I\sqcup J=\{1,\dots,n\}}
\int_{\mgnbar{q}{h+1}}\prod_{i\in I}\psi_i^{d_i}
\int_{\mgnbar{r}{k+1}}\kappa_1\prod_{j\in J}\psi_j^{d_j}\cr
&=\sum_{0\leq q\leq g\atop I\subset \{1,\dots,n\}}
\int_{\mgnbar{g}{n}}{\xi_{q,I}}_*(1\times\kappa_1)\prod_{i=1}^n\psi_i^{d_i}
\,.}
$$

\medskip\noindent{- Term 5.}
$$
\langle\tau_0\tau_1\tau_{\ul{d}}\rangle=
\langle\tau_0\tau_1\tau_{\ul{d}}\rangle_{g-1,n+2}=
\int_{\mgnbar{g}{n}}{\xi_{irr}}_*(\psi_{n+2})\prod_{i=1}^n\psi_i^{d_i}\,.
$$
This is an expression which appeared also in the formula for
$\wmn{(0,m_1,0,1,0,\dots)}{n}$. As in that case, it could have been
interpreted,
alternatively, as
$$
24\int_{\mgnbar{g}{n}}{\xi_{1,\emptyset}}_*(\psi_1\times\kappa_0)
\prod_{i=1}^n\psi_i^{d_i}\,.
$$

\medskip\noindent{- Term 6.} The relevant part is
$$
\sum_{I\sqcup J=\{1,\dots,n\}}\hskip-0.4truecm
\langle\tau_0\tau_{\ul{d}_I}\rangle_{q,h+1}
\langle\tau_1\tau_{\ul{d}_J}\rangle_{r,k+1}\,,
$$
where $q+r=g$, and thus it can be rewritten as
$$
\sum_{0\leq q\leq g\atop I\subset \{1,\dots,n\}}
\int_{\mgnbar{g}{n}}{\xi_{q,I}}_*(1\times\psi_{k+1})
\prod_{i=1}^n\psi_i^{d_i}\,.
$$

\medskip\noindent{- Term 7.} We have
$$
\langle\tau_0^4\tau_{\ul{d}}\rangle=
\langle\tau_0^4\tau_{\ul{d}}\rangle_{g-2,n+4}=
\int_{\mgnbar{g}{n}}{\xi_A}_*(1)\prod_{i=1}^n\psi_i^{d_i}\,,
$$
where the graph $A$ is depicted in Figure 2.

\medskip\noindent{- Term 8.} This is
$$
\sum_{I\sqcup J=\{1,\dots,n\}}\hskip-0.4truecm
\langle\tau_0^3\tau_{\ul{d}_I}\rangle_{q,h+3}
\langle\tau_0\tau_{\ul{d}_J}\rangle_{r,k+1}\,,
$$
with $q+r=g-1$. Clearly it can also be written as
$$
\sum_{p,P}\int_{\mgnbar{g}{n}}{\xi_{(B,p,P)}}_*(1)\prod_{i=1}^n\psi_i^{d_i}\,,
$$
where the directed graph $B$ is illustrated in Figure 2, and $p$ (resp., $P$)
runs
through all possible assignments of genera to the vertices of $B$ subject to
the
condition that their sum be equal to $g-1$ (resp., all partitions of
$\{1,\dots,n\}$ indexed by the vertices of $B$).

\medskip\noindent{- Term 9.} This is
$$
\sum_{I\sqcup J\sqcup K=\{1,\dots,n\}}
\hskip-0.4truecm
\langle\tau_0^2\tau_{\ul{d}_I}\rangle_{q,h+2}
\langle\tau_0\tau_{\ul{d}_J}\rangle_{r,k+1}
\langle\tau_0\tau_{\ul{d}_K}\rangle_{s,l+1}\,,
$$
where $q$, $r$ and $s$ add to $g$. As in the preceding case, then, this term
can
also be written
$$
\sum_{p,P}\int_{\mgnbar{g}{n}}{\xi_{(C,p,P)}}_*(1)\prod_{i=1}^n\psi_i^{d_i}\,,
$$
where the directed graph $C$ is illustrated in Figure 2, and $p$ (resp., $P$)
runs through all possible assignments of genera to the vertices of $C$ subject
to
the condition that their sum be equal to $g$ (resp., all partitions of
$\{1,\dots,n\}$ indexed by the vertices of $C$).

\medskip\noindent{- Term 10.} This is handled similarly to the two preceding
ones. Consider the directed graph $D$ in Figure 2. Then
$$
\sum_{I\sqcup J=\{1,\dots,n\}}
\hskip-0.4truecm
\langle\tau_0^2\tau_{\ul{d}_I}\rangle
\langle\tau_0^2\tau_{\ul{d}_J}\rangle
=\sum_{p,P}\int_{\mgnbar{g}{n}}{\xi_{(D,p,P)}}_*(1)\prod_{i=1}^n\psi_i^{d_i}\,,
$$
where $p$ runs through all assignments of genera to the vertices of $D$ adding
to $g-1$ and $P$ through all partitions of $\{1,\dots,n\}$ indexed by the
vertices
of $D$.

\medskip\noindent{- Term 11.} The expression
$\langle\tau_0\tau_{\ul{d}}\rangle$ already appears in in the formula for
$\wmn{(0,m_1,0,1,0,\dots)}{n}$, and we have seen that it equals
$$
24\int_{\mgnbar{g}{n}}{\xi_{1,\emptyset}}_*(\psi_1\times 1)
\prod_{i=1}^n\psi_i^{d_i}\,.
$$
What all the above computation suggests is
that a reasonable candidate for an expression of $\wmn{(0,m_1,2,0,\dots)}{n}$
in
terms of the standard algebro-geometric classes might be
$$
\eqalign{&\quad 72\kappa_1^2-348\kappa_2
-6{\xi_{irr}}_*(\kappa_1)-12\hskip-0.3truecm
\sum_{0\leq q\leq g\atop I\subset \{1,\dots,n\}}
\hskip-0.3truecm{\xi_{q,I}}_*(1\times\kappa_1)
+12{\xi_{irr}}_*(\psi_{n+2})\cr
&+12\hskip-0.3truecm\sum_{0\leq q\leq g\atop I\subset \{1,\dots,n\}}
\hskip-0.3truecm{\xi_{q,I}}_*(1\times\psi_{k+1}) +{1\over
8}{\xi_A}_*(1)+{1\over 2}\sum_{p,P}{\xi_{(B,p,P)}}_*(1)\cr
&+{1\over 2}\sum_{p,P}{\xi_{(C,p,P)}}_*(1)+{1\over
4}\sum_{p,P}{\xi_{(D,p,P)}}_*(1) -36{\xi_{1,\emptyset}}_*(\psi_1\times 1)\,,
}\leqno{(3.15)}$$
up to the ambiguity noticed in the analysis of term 5.

As we remarked after formula (3.11), due to a number of remarkable
cancellations, in the expressions of the derivatives of $F$ with respect to the
$s$
variables in terms of derivatives with respect to the $t$ variables many
derivatives of $F$ that would a priori be allowed by the DFIZ theorem are not
present. For example, in codimension 2 the terms
$$
\left(\partial F\over\partial t_0\right)^4\,,\quad
\left(\partial F\over\partial t_2\right)^2\,,
$$
among many others, are missing. This is wonderful, for otherwise our formulas
would have been ruined. In fact, when studying codimension 2
classes in genus $g$, we would have gotten terms of the sort
$$
\eqalign{
&\langle\tau_0\tau_{\ul{d}_{I_1}}\rangle_{q_1,h_1+1}
\langle\tau_0\tau_{\ul{d}_{I_2}}\rangle_{q_2,h_2+1}
\langle\tau_0\tau_{\ul{d}_{I_3}}\rangle_{q_3,h_3+1}
\langle\tau_0\tau_{\ul{d}_{I_4}}\rangle_{q_4,h_4+1}\,,\cr
&\langle\tau_2\tau_{\ul{d}_{I_1}}\rangle_{r_1,h_1+1}
\langle\tau_2\tau_{\ul{d}_{I_2}}\rangle_{r_2,h_2+1}\,,
}$$
respectively. Now it is easy to see that we must have
$q_1+q_2+q_3+q_4=g+1=r_1+r_2$. Thus it would have been impossible to write
these terms as intersection numbers on $\mgnbar{g}{n}$.

In the same vein, but with considerably more effort, we could have given
similar formulas for some classes $\wmn{m_*}{n}$ of higher codimension,
including in particular all those of codimension 3. In the Appendix we have
listed the expressions of the derivatives of $F$ with respect to the $s$
variables
in terms of those with respect to the $t$ variables that are needed to carry
out
these computations, which are otherwise left to the reader.

The resulting identities are of the form
$$
\int_{\wmn{m_*}{n}}\prod\psi_i^{d_i}
=\int_{\mgnbar{g}{n}}X_{m_*,n}\prod\psi_i^{d_i}\,,
$$
where $g$ is given by $4g-4+2n=\sum m_i(2i-1)$ and $X_{m_*,n}$ is a polynomial
in the
Mumford classes and the boundary classes. The reader is invited to check that
this is
indeed true, using the methods employed in this section. Of
course, the point is to verify that all the terms that one gets can be
interpreted as
intersection numbers on $\mgnbar{g}{n}$, and not in higher genus. The formulas
one
finds, arranged by increasing codimension, look as follows
\bigskip
\noindent - Codimension 1
$$
X_{(0,m_1,1,0,\dots),n}=12\kappa_1+\cdots
$$
- Codimension 2
$$
\eqalign{
&X_{(0,m_1,0,1,0,\dots),n}=120\kappa_2+\cdots\cr
&X_{(0,m_1,2,0,\dots),n}=72\kappa_1^2-348\kappa_2+\cdots
}$$
- Codimension 3
$$
\eqalign{
&X_{(0,m_1,0,0,1,0,\dots),n}=1680\kappa_3+\cdots\cr
&X_{(0,m_1,1,1,0,\dots),n}=1440\kappa_1\kappa_2-13680\kappa_3+\cdots\cr
&X_{(0,m_1,3,0,\dots),n}=288\kappa_1^3-4176\kappa_1\kappa_2
+20736\kappa_3+\cdots
}$$
- Codimension 4
$$
\eqalign{
&X_{(0,m_1,0,0,0,1,0,\dots),n}=30240\kappa_4+\cdots\cr
&X_{(0,m_1,1,0,1,0\dots),n}=20160\kappa_1\kappa_3-312480\kappa_4+\cdots\cr
&X_{(0,m_1,0,2,0\dots),n}=7200\kappa_2^2-159120\kappa_4+\cdots
}$$
- Codimension 5
$$
X_{(0,m_1,0,0,0,0,1,0,\dots),n}=665280\kappa_5+\cdots
$$
- Codimension 6
$$
X_{(0,m_1,0,0,0,0,0,1,0,\dots),n}=17297280\kappa_6+\cdots
$$
The dots stand for boundary classes, which are well determined up to
ambiguities of
the kinds previously described. This list is complete up to codimension 3
included.

The same remarkable cancellations that we observed for codimension two classes
occur, to an even greater extent, in higher codimension. For instance, the
expression
for $\partial^3 F/\partial s_2^3$ given in the appendix involves 41 terms
while, a
priori, up to 585 might have been expected from the statement of the DFIZ
theorem.

\vskip1.5 truecm
\specialhead 4. The codimension one case
\endspecialhead
Our main goal in this section is to complete the study of the codimension
one class $\wmn{(0,m_1,1,0,\dots)}{n}$. For simplicity, this will be
denoted simply by $W$ throughout the section. We shall prove the following

\proclaim{Proposition 1}When $n\geq 2$, for any class $\gamma\in
H^{6g-6+2n-2}(\mprime{g}{n},\Bbb Q)$ one has
$$
\int_{W}\gamma=\int_{\mgnbar{g}{n}}\alpha^*(\gamma)
(12\kappa_1-\delta)\,,
$$
where $\alpha$ is the natural map from $\mgnbar{g}{n}$ to $\mprime{g}{n}$.
\endproclaim

In section 3 we have shown that the proposition holds when $\gamma$
is a product of classes $\psi_i$. Notice that both
$W$ and
$12\kappa_1-\delta$ are invariant under the natural action of the symmetric
group $\Cal S_n$ on
$\mgnbar{g}{n}$. Moreover, it is reasonable to expect that $W$ can be
``lifted'',
non-uniquely, to a homology class on $\mgnbar{g}{n}$. Proving this, however,
requires a little argument that will be given later. Granting this, the
proposition is then a direct consequence of the following lemma.

\proclaim{Lemma 2}Let $x\in H^2(\mgnbar{g}{n},\Bbb Q)$ be an $\Cal
S_n$ invariant class, with $n\geq 2$. Suppose that
$$
\int_{\mgnbar{g}{n}}x\cup\prod_{i=1}^n
\psi_i^{d_i}=0
$$
for all choices of $d_1,\dots,d_n$. Then $x$ is a linear combination of classes
$\delta_{p,\emptyset}$.
\endproclaim

Incidentally, it is obvious that the converse of the statement of the lemma is
true. Moreover, the $\delta_{p,\emptyset}$ are precisely the classes of those
components of the boundary that are partially contracted by $\alpha$. We
now prove the lemma.

It follows from a well-known theorem of Harer [4] that the second
cohomology group of $\mgnbar{g}{n}$ is generated by the
classes $\kappa_1,\psi_1,\dots,\psi_n$ and by the boundary classes
$\delta_{\Gamma}$, where $\Gamma$ runs through all isomorphisms
classes of dual graphs having only one edge. We set
$\psi=\sum_{i=1}^n\psi_i$. Given integers $p$ and $h$, with
$0\leq p\leq g$ and $0\leq h\leq n$, we denote by $\delta_{p,h}$ the sum
$\sum\delta_\Gamma$, where $\Gamma$ runs through all isomorphism
classes of dual graphs of curves with exactly two components meeting at one
point, one of which has genus $p$ and carries $h$ marked points. In terms of
these, one may write the class of the boundary as
$$
\delta=\hskip-0.4truecm\sum_{0\leq p\leq g/2\atop{0\leq h\leq n\atop h\leq
n/2\text{ if }p=g/2}}\hskip-0.4truecm
\delta_{p,h}\,.
$$
Clearly, any invariant class in the second cohomology group of
$\mgnbar{g}{n}$ is a linear combination of $\kappa_1$, $\psi$, $\delta_{irr}$,
and the $\delta_{p,h}$. If $g=2$ (resp., $g=1$, resp., $g=0$) one can do
without
$\kappa_1$ (resp., $\kappa_1$ and $\psi$, resp., $\kappa_1$, $\psi$ and
$\delta_{irr}$). We then write
$$
x=a\kappa_1+b\psi+c\delta_{irr}+\sum c_{p,h}\delta_{p,h}\,.
$$
We will first show that $a=b=c=0$. Set
$$
\alpha_s=\prod_{i=1}^{n-2}\psi_i\psi_{n-1}^s\psi_n^d\,,
$$
where $d=3g-s-2$, and notice that, if $s\not\equiv 2\mod{3}$, then
$\int\delta_{p,h}\alpha_s=0$. In fact, $\int\delta_{p,h}\alpha_s$ is a linear
combination of terms of the form
$\ang{\tau_0\tau_1^a}\ang{\tau_0\tau_1^b\tau_s\tau_d}$ or of the form
$\ang{\tau_0\tau_1^a\tau_s}\ang{\tau_0\tau_1^b\tau_d}$. But since $a-(a+1)$
and $a+s-(a+2)$ are not divisible by 3 ($s\not\equiv 2\mod{3}$), both
$\ang{\tau_0\tau_1^a}$ and $\ang{\tau_0\tau_1^a\tau_s}$ vanish. We now
wish to compute $\int\kappa_1\alpha_s$, $\int\psi\alpha_s$, and
$\int\delta_{irr}\alpha_s$. For this we are going to use the following
well-known formulae [16].
$$
\eqalign{
\ang{\tau_0\tau_{d_1}\dots\tau_{d_n}}&=\sum_{d_i>0}
\ang{\tau_{d_1}\dots\tau_{d_i-1}\dots\tau_{d_n}}\,,\cr
\ang{\tau_1\tau_{d_1}\dots\tau_{d_n}}&=(2g-2+n)
\ang{\tau_{d_1}\dots\tau_{d_n}}\,,
}\leqno{(4.1)}$$
where $3g-3+n=\sum d_i$. The first formula is just the
string equation (1.9), while the second is a special case of (1.7). Setting
$r=2g-2$ we then have
$$
\int\kappa_1\alpha_s=\ang{\tau_2\tau_1^{n-2}\tau_s\tau_d}
={(r+n)!\over (r+2)!}\ang{\tau_2\tau_s\tau_d}_{g,3}\,,
$$
$$
\eqalign{
\int\psi\alpha_s&=(n-2)\ang{\tau_2\tau_1^{n-3}\tau_s\tau_d}+
\ang{\tau_1^{n-2}\tau_{s+1}\tau_d}+\ang{\tau_1^{n-2}\tau_s\tau_{d+1}}\cr
&={(r+n-1)!\over (r+2)!}(n-2)\ang{\tau_2\tau_s\tau_d}_{g,3}+{(r+n-1)!\over
(r+1)!}\ang{\tau_{s+1}\tau_d}_{g,2}\cr
&\quad+{(r+n-1)!\over (r+1)!}\ang{\tau_s\tau_{d+1}}_{g,2}\,,
}$$
$$
\eqalign{
2\int\delta_{irr}\alpha_s&=\ang{\tau_0^2\tau_1^{n-2}\tau_s\tau_d}
={(r+n-1)!\over (r+1)!}\ang{\tau_0^2\tau_s\tau_d}\cr
&={(r+n-1)!\over(r+1)!}\left(\ang{\tau_{s-2}\tau_d}_{g-1,2}+
2\ang{\tau_{s-1}\tau_{d-1}}_{g-1,2}+
\ang{\tau_s\tau_{d-2}}_{g-1,2}\right)\,.
}$$
To simplify these expressions we shall use the fundamental fact [7] that
the function $Z(t_*)=\exp F(t_*)$ satisfies the KdV equation. It is convenient
to
set
$$
\eqalign{
&\varphi(g)=\left\{\eqalign{
&\ang{\tau_{3g-2}}\qquad\text{if }g\geq 1\,,\cr
&\ang{\tau_0^3}=1\qquad\,\,\,\,\text{if }g=0\,,
}\right.\cr
&\phantom{.}\cr
&T(s,r)=\frac{\ang{\tau_s\tau_d}_{g,2}}{\varphi(g)}\,,\cr
&\phantom{.}\cr
&U(s,r)=\frac{\ang{\tau_2\tau_s\tau_d}_{g,3}}{\varphi(g)}\,.
}$$
Writing the KdV equation in Lax
form (also referred to as Gel'fand-Dikii form), one gets in particular (cf.
[16], page 251) that
$$
\varphi(g)=\frac{1}{24g}\varphi(g-1)=\frac{1}{12(r+2)}\varphi(g-1)\,.
$$
It follows that
$$
\left\{\eqalign{
\int\kappa_1\alpha_s&={(r+n-1)!\over (r+2)!}\varphi(g)(r+n)U(s,r)\,,\cr
\int\psi\alpha_s&={(r+n-1)!\over (r+2)!}\varphi(g)
\left[(n-2)U(s,r)\right.\cr
&\quad\left.+(r+2)T(s+1,r)+(r+2)T(s,r)\right]\,,\cr
\int\delta_{irr}\alpha_s&={(r+n-1)!\over (r+2)!}\varphi(g)6(r+2)^2\left[
T(s-2,r-2)\right.\cr
&\quad\left.+2T(s-1,r-2)+T(s,r-2)\right]\,.
}\right.\leqno{(4.2)}
$$
To calculate $T(s,r)$ and $U(s,r)$ we use again the fact that $Z$ satisfies the
KdV equation, but expressing this by saying that $Z$ is annihilated by the
Virasoro operators $L_k$, for $k\geq -1$. The equations $L_{-1}Z=L_0Z=0$
are equivalent to the (4.1). The Virasoro operator $L_k$ is given, for $k>0$,
by
$$
\eqalign{
L_k&=-\frac{(2k+3)!!}{2}{\partial\over\partial t_{k+1}}+
\frac{1}{2}\sum_{i=0}^{\infty}(2k+2i+1)(2k+2i-1)\dots(2i+1)t_i
\dfrac{\partial}{\partial t_{i+k}}\cr
&\quad+\frac{1}{4}\sum_{r+s+1=k}(2r+1)!!(2s+1)!!{\partial^2\over\partial
t_r\partial t_s}
}$$
Recalling that
$$
F(t_*)=\sum\ang{\tau_{d_1}\dots\tau_{d_n}}t_{d_1}\dots t_{d_n}\,,
$$
to say that $L_kZ=0$ translates into
$$
\eqalign{
\ang{\tau_{k+1}\tau_{\ul{d}}}&=\frac{1}{(2k+3)!!}\left[\sum_{j=1}^n
\frac{(2k+2d_j+1)!!}{(2d_j-1)!!}\ang{\tau_{d_1}\dots
\tau_{d_{j+k}}\dots\tau_{d_n}}\right.\cr
&\quad+\frac{1}{2}\sum_{r+s=k-1}
(2r+1)!!(2s+1)!!\ang{\tau_r\tau_s\tau_{\ul{d}}}\cr
&\left.\quad+\frac{1}{2}\sum_{r+s=k-1} (2r+1)!!(2s+1)!!
\sum_{I\subset\{1,\dots,n\}}\ang{\tau_r\tau_{\ul{d}_I}}
\ang{\tau_s\tau_{\ul{d}_{\Cal CI}}}\right]
}$$
for any $\ul{d}=(d_1,\dots,d_n)$. It follows that, provided $s\not\equiv
2\mod 3$,
$$
\eqalign{
U(s,r)&=\frac{1}{3\cdot
5}\left[(2s+3)(2s+1)T(s+1,r)+(3r-2s+5)(3r-2s+3)T(s,r)\right.\cr
&\quad\left. +6(r+2)(T(s-2,r-2)+2T(s-1,r-2)+T(s,r-2))\right]\,.
}$$
Furthermore we have
$$
\eqalign{
T(s,r)&=0\qquad\text{if }s<0\,, \cr
T(0,r)&=1\,, \cr T(1,r)&=r+1\,, \cr
T(2,r)&=\frac{1}{3\cdot 5}\left[(3r+3)(3r+1) +6(r+2)\right]\,,\cr
T(3,r)&=\frac{1}{3\cdot 5\cdot 7}\left[(3r+3)(3r+1)(3r-1) +3\cdot 12
r(r+2)+\frac{3}{2}(r+2)\right]\,,\cr
T(4,r)&=\frac{1}{3\cdot 5\cdot 7\cdot
9}\left[(3r+3)(3r+1)(3r-1)(3r-3)\phantom{\frac{1}{1}}\right.\cr
&\quad\left.  +12(r+2)\left(3\cdot 5+\frac{9}{2}r+\frac{3}{8}
\right)T(1,r-2)+3\cdot 5\cdot 12(r+2)T(2,r-2)\right]\,,\cr
T(5,r)&=\frac{1}{3\cdot 5\cdot 7\cdot 9\cdot
11}\left[(3r+3)(3r+1)(3r-1)(3r-3)(3r-5)\phantom{\frac{1}{1}}\right.\cr
&\quad\left. +12(r+2)\left(\left(3\cdot 5\cdot 7+3\cdot 3\cdot 5\cdot
r+\frac{3\cdot 5}{8}\right)T(2,r-2)+3\cdot 5\cdot
7\cdot T(3,r-2)\right)\right]\,.}
$$
Now consider the system of three linear equations in the unknowns $a$, $b$
and $c$ given by
$$
\frac{(r+2)!}{\varphi(g)(r+n-1)!}\int x\alpha_s=0\,,\qquad s=0,1,3\,.
\leqno{\text{Eq}_s)}
$$
The coefficients are implicitly given by (4.2), an can be calculated using the
formulas for $U(s,r)$ and $T(s,r)$ we have just given. An algebraic
calculation (best done by computer) shows that the determinant of this
system equals
$$
\frac{36}{875}(r-2)r^2(r+2)^5(n+r)(4r+17)\,.
$$
Recall that $n$ is an integer greater or equal to 2, and that $r$ is an integer
greater or equal to $-2$; thus our determinant vanishes only for $r=-2,0,2$,
that is, for $g=0,1,2$. This shows that $a=b=c=0$ for $g\geq 3$. If $g=2$, the
class $\kappa_1$ is linearly dependent on the others, so that we may set
$a=0$ and view $\text{Eq}_0$ and $\text{Eq}_1$ as a linear system in the
unknowns $b$ and $c$. The determinant of this system equals
$(1536/5)(n+2)$, which is non-zero in our situation. If $g=1$ one may set
$a=b=0$, and the coefficient of $c$ in $\text{Eq}_0$ is $24$. If $g=0$ the
classes $\kappa_1$, $\psi$ and $\delta_{irr}\ (=0)$ are linear combinations of
the $\delta_{p,h}$.

We have thus shown that, for any value of the genus $g$, the class $x$ is a
linear combination
$$
x=\sum c_{p,h}\delta_{p,h}\,.
$$
We wish to show that $c_{p,h}$ vanishes unless $h=0$ or $h=n$. Set
$$
\eqalign{
\beta_{0,j}&=\psi_1^{s-j}\prod_{h=2}^{j-1}\psi_h\,,\qquad s=3g-2+n\,,\cr
\beta_{q,j}&=\psi_1^{s-j}\psi_2^{3q-1}\prod_{h=3}^{j+1}\psi_h\,,\qquad
s=3(g-q)-2+n\qquad\text{if }q>0\,.}
$$
In what follows we shall always assume, as we may, that $p\leq g/2$. The
intersection number $\int \delta_{p,h}\beta_{0,j}$ is, a priori, a linear
combinations of terms of the form
$$
\ang{\tau_0^{a+1}\tau_1^c}_{p,h+1}
\ang{\tau_0^{b+1}\tau_1^d\tau_{3g-2+n-j}}_{g-p,n-h+1}
$$
or
$$
\ang{\tau_0^{a+1}\tau_1^c\tau_{3g-2+n-j}}_{p,h+1}
\ang{\tau_0^{b+1}\tau_1^d}_{g-p,n-h+1}\,.
$$
However, since $\ang{\tau_0^{b+1}\tau_1^d}_{g-p,n-h+1}$ is non-zero only if
$b=2$, in which case $g-p=0$, there are no terms of the second kind. As for
those of the first kind, they may be non-zero only if $a=2$, $p=0$, $h=a+c$,
$c+d=j-2$, so that $h\leq a+j-2=j$. In conclusion
$$
\int \delta_{p,h}\beta_{0,j}=0\qquad \text{if }p>0\text{ or }p=0,\,h>j\,;
\leqno{(4.3)}
$$
moreover
$$
\int \delta_{0,j}\beta_{0,j}\neq 0\qquad \text{if }1<j<n\,,\leqno{(4.4)}
$$
since this number is a multiple, with positive coefficient, of
$\ang{\tau_0^3\tau_1^{j-2}}
\ang{\tau_0^{n-j}\tau_{3g-2+n-j}}\neq 0$.
Let us now compute the intersection number
$$
\int \delta_{p,h}\beta_{q,j}
$$
when $q>0$ and, as usual, $p\leq g/2$. This is, a priori, a linear combination
of terms of the form
$$
\ang{\tau_0^{a+1}\tau_1^c\tau_{3q-1}}_{p,h+1}
\ang{\tau_0^{b+1}\tau_1^d\tau_{s-j}}_{g-p,n-h+1}\,,
$$
or
$$
\ang{\tau_0^{a+1}\tau_1^c\tau_{s-j}}_{p,h+1}
\ang{\tau_0^{b+1}\tau_1^d\tau_{3q-1}}_{g-p,n-h+1}\,,
$$
or else
$$
\ang{\tau_0^{a+1}\tau_1^c}_{p,h+1}
\ang{\tau_0^{b+1}\tau_1^d\tau_{3q-1}\tau_{s-j}}_{g-p,n-h+1}\,,
$$
or, finally,
$$
\ang{\tau_0^{a+1}\tau_1^c\tau_{3q-1}\tau_{s-j}}_{p,h+1}
\ang{\tau_0^{b+1}\tau_1^d}_{g-p,n-h+1}\,.
$$
We already saw that a term of this last type is not equal to zero only
if $g-p=0$ , which is impossible. Terms of the third
type are different from zero only if $p=0$, $a=2$, $h=c+2\leq j+1$. Let us
analyze terms of the first  type. These are non-zero only if $3p-3+h+1=3q-1-c$.
On the other hand $h=a+c+1$, so that $3q=3p+a$. Furthermore $c+d=j-1$, so
that $h\leq j+a$. The same argument shows that there are no non-zero terms
of the second type. We conclude that
$$
\int \delta_{p,h}\beta_{q,j}=0\qquad \text{if }p>q\text{ or }p\leq q,
\,h>j+3(q-p)\,;\leqno{(4.5)}
$$
moreover
$$
\int \delta_{p,h}\beta_{p,h}\neq 0\qquad \text{if }0<h<n\,,\leqno{(4.6)}
$$
since this number is a positive multiple of
$\ang{\tau_0\tau_1^{h-1}\tau_{3p-1}}
\ang{\tau_0^{n-h}\tau_{3(g-p)-3+n-h+1}}\neq 0$. Arguing by double
induction on $p$ and $h$, it follows from (4.3), (4.4), (4.5), and (4.6) that
$c_{p,h}=0$ for all $p$ and all $h$ different from 0 and $n$. This proves the
lemma.

\bigskip
As we announced, to complete the proof of Proposition 1 it remains to
compare the (co)homology of $\mprime{g}{n}$ with that of $\mgnbar{g}{n}$.
Rational coefficients will be used throughout. We shall show that

\proclaim{Lemma 3}There is an exact sequence
$$
0\to H^{6g-6+2n-2}(\mprime{g}{n})\mathop{\longrightarrow}^{\alpha^*}
H^{6g-6+2n-2}(\mgnbar{g}{n})\to A\to 0\,,
$$
where $\alpha:\mgnbar{g}{n}\to\mprime{g}{n}$ is the natural map and $A$
is the vector space freely generated by the boundary classes
$\delta_{p,\emptyset}$, with $1\leq p\leq g$ \rm{(}or $1\leq p\leq g-1$ for
$n=1$\rm{)}.
\endproclaim

A consequence of this lemma is that the functional defined by integration on
$W$ lifts to an element $\widetilde{W}$ of
$H^{6g-6+2n-2}(\mgnbar{g}{n})^\vee\cong H_{6g-6+2n-2}(\mgnbar{g}{n})$,
which we may choose to be $\Cal S_n$-invariant, since $W$ is. Proposition 1
follows by applying Lemma 2 to the difference between $12\kappa_1-\delta$
and the Poincar\'e dual of $\widetilde{W}$.

We now prove the lemma. Look at the commutative diagram
$$
\CD H_{\text{c}}^{d-2}(\mp\setminus\Sigma')  @>\rho>>  H^{d-2}(\mp)
@>{\alpha^*}>>  H^{d-2}(\mm)\\ @|  @.  @|\\
H_{\text{c}}^{d-2}(\mp\setminus\Sigma')  @<{\cong}<<
H_{\text{c}}^{d-2}(\mm\setminus\Sigma)  @>>> H^{d-2}(\mm)\\ @VVV  @VVV
@VVV\\ H_2(\mp\setminus\Sigma')  @<<< H_2(\mm\setminus\Sigma)  @>\sigma>>
H_2(\mm)
\endCD
$$
where we have set $d=6g-6+2n$, $\mm=\mgnbar{g}{n}$,
$\mp=\mprime{g}{n}$, and $\Sigma$ (resp., $\Sigma'$) stands for the union of
the components of $\partial\mm$ of the form $\Delta_{p,\emptyset}$ (resp.,
the image of $\Sigma$ in $\mp$). The three vertical arrows are isomorphisms
by Poincar\'e duality. The map $\rho$ is a piece of the exact sequence of
cohomology with compact support
$$
\cdots\to H^{d-3}(\Sigma')\to
H_{\text{c}}^{d-2}(\mp\setminus\Sigma')\mathop{\longrightarrow}^{\rho}
H^{d-2}(\mp)\to H^{d-2}(\Sigma')\to\cdots\,.
$$
On the other hand $H^{d-3}(\Sigma')$ and $H^{d-2}(\Sigma')$ vanish since
the dimension of $\Sigma'$ is strictly smaller than $d-3$, so $\rho$ is an
isomorphism. It follows in particular that $\alpha^*$ is injective if and only
if
$\sigma$ is, and that its cokernel can be identified with the one of $\sigma$.
Passing to duals, we have to look at
$\sigma^\vee : H^2(\mm)\to H^2(\mm\setminus\Sigma)$, which fits into the
exact sequence
$$
\cdots\to
H^2(\mm;\mm\setminus\Sigma)\to
H^2(\mm)\mathop{\longrightarrow}^{\sigma^\vee}
H^2(\mm\setminus\Sigma)\to
H^3(\mm;\mm\setminus\Sigma)\to\cdots\,.
$$
Now the Thom isomorphism implies that $H^3(\mm;\mm\setminus\Sigma)$
vanishes, since $H^1(\mgnbar{p}{h})$ does, for any $p$ and $h$, and that,
moreover, $H^2(\mm;\mm\setminus\Sigma)$ is freely generated by the
classes of the components of $\Sigma$. Since the images of these are
independent in $H^2(\mm)$, the conclusion follows.

\vskip1.5 truecm
\specialhead 5. Examples and comments
\endspecialhead
It is instructive to work out a couple of simple examples. It should be
clear from them how intricate a direct attack on the problem would be. We
begin by checking formula (2.3) on $\mgnbar{1}{1}$ and on $\mgnbar{0}{4}$. One
thing that simplifies matters in these cases is that, for these values of $g$
and $n$,
one has $\mprime{g}{n}=\mgnbar{g}{n}$. We shall write $M$ for
$\mgn{(0,m_1,1,0,\dots)}{n}$ and $W$ for $\wmn{(0,m_1,1,0,\dots)}{n}$. In
general, if $v$ and $l$ stand for the numbers of vertices and edges of a ribbon
graph of genus $g$ with $n$ boundary components all of whose vertices are
trivalent save for a pentavalent one, one has
$$ v=2n+4g-6\,,\qquad l=3n+6g-8\,.
$$ In the case of $\mgnbar{1}{1}$ this yields $v=0$; thus, $W=0$. To check
(2.3) in
this case we must then show that $12\kappa_1=\delta$. But now one knows that
$\tilde{\kappa}_1=12\lambda-\delta$ vanishes on $\mgnbar{1}{1}$
(cf. [6], for instance). On the other hand $\psi=\lambda$. Thus
$12\kappa_1=12\psi=12\lambda=\delta$, as desired.

The case of $\mgnbar{0}{4}$ is more entertaining. The formulas above give
$v=2$,
$l=4$; in particular, $M$ is 4-dimensional and hence $W$ zero-dimensional. The
possible graphs in $M$ are those of type a), b), and c) in Figure 3; their
degenerations are illustrated in d), e), f), and g). Among these, the first two
are
internal to moduli, while the last two correspond to points in the boundary.
Now let us consider the projection
$$
\eta: \mcombar{0}{4}=\mgnbar{0}{4}\times\Bbb R^4\to\Bbb R^4
$$ which associates to any numbered ribbon graph with metric the quadruple
of positive real numbers given by the lengths of its four
boundary components. Clearly, a cycle $Z$ in $\mgnbar{0}{4}$ representing $W$
can be obtained by cutting $M$ with a section $\eta=(P_1,\dots,P_4)$, where the
$P_i$ are positive constants.
We choose $P_1,\dots,P_4$ in such a way that $P_{i+1}\geq 10P_i$, for
$i=1,2,3$.
Since for graphs f) and g) two of the perimeters necessarily coincide,
Z is entirely contained in the interior of moduli. We now show that graphs of
types
c), d), and e) cannot occur in $Z$. In fact, for graphs of type d) one of the
perimeters equals the sum of the remaining three, and this is forbidden by our
choice of $P$'s. In e) one of the perimeters, which we may assume to be the
longest, equals the sum of two of the other perimeters minus the remaining one.
This too is incompatible with our choices. Let us now look at graph c), where
the
edges have been labelled with their respective lengths $l_1,\dots,l_4$. Up to
the
numbering of the boundary components the perimeters are
$$ p_1=l_1+l_4\,,\quad p_2=l_2+l_4\,,\quad p_3=l_3\,,\quad p_4=l_1+l_2+l_3\,.
$$ We also have the obvious inequalities
$$ p_4<p_1+p_2+p_3\,,\quad p_3<p_4\,,\quad p_1<p_2+p_4\,,\quad p_2<p_1+p_4\,.
$$ These inequalities imply, in order, that neither $p_4$, nor $p_3$, nor $p_2$
or
$p_1$ can equal $P_4$. This excludes case c). We next examine case a). It is
clear
that the longest perimeter is the ``external'' one and that, except for this
restriction, the perimeters can be arbitrarily assigned. Therefore this case
accounts for $6=3!$ points of $Z$, one for each choice of labelling of the
three
``internal'' boundary components by
$\{1,2,3\}$. We now examine graph b). Up to the numbering of the boundary
components the perimeters are
$$
p_1=l_1\,,\quad p_2=l_2\,,\quad p_3=l_3+l_1\,,\quad p_4=l_2+l_3+2l_4\,.
$$
The following inequalities hold
$$
p_1<p_3\,,\quad p_2<p_4\,,\quad p_3<p_1+p_4\,,\quad p_1+p_2<p_3+p_4\,.
$$
{}From the first three inequalities we get in particular that we must have
$p_4=P_4$. The only possibilities for $(p_1,p_2,p_3,p_4)$ are
$$
(P_1,P_2,P_3,P_4)\,,\quad (P_1,P_3,P_2,P_4)\,,\quad (P_2,P_1,P_3,P_4)\,.
$$
In conclusion, the support of $Z$ consists of 9 points. We claim that $Z$ is
the sum of
these points, taken with the positive sign. This follows immediately from
Kontsevich's
recipe (cf. [7], page 11) for the orientation of $Z$. In his notation, this is
given by
$\Omega^d$, where $d$ is the dimension of $Z$, i.e., by the constant 1.
It follows that $W$ is 9 times the
fundamental class of
$\mgnbar{0}{4}=\Bbb P^1$. To prove formula (2.3) in the present case it now
suffices to show that $12\kappa_1-\delta$ has degree 9. This follows if we can
show that
$$
\deg\delta=3\,,\quad\deg\psi=4\,,\quad\deg\tilde{\kappa}_1=-\deg\delta=-3\,.
$$ We briefly indicate how to do it. Set $X'=\Bbb P^1\times\Bbb P^1$, denote by
$f'$ the projection onto the second factor, and by $\sigma'_i$, $i=1,2,3$,
three
constant sections of $f'$. Blow up $X'$ at the three points where these
sections
meet the diagonal, to obtain a family $f:X\to\Bbb P^1$ together with four
distinct
sections
$\sigma_1,\dots,\sigma_4$, which are the proper transforms of
$\sigma'_1,\dots,\sigma'_3$ and of the diagonal. This is the universal curve
over
$\mgnbar{0}{4}$. Denote by $E_1,E_2,E_3$ the exceptional curves of the blow-up,
and set $D_i=\sigma_i(\Bbb P^1)$. Then
$$
\eqalign{ 0&=(D_i+E_i)^2=D_i^2+2+E_i^2=D_i^2+1\qquad\text{if }i\leq 3\cr
2&=(D_4+\textstyle{\sum}E_i)^2=D_4^2+6-3=D_4^2+3\,. }
$$
In all cases $\deg\psi_i=-D_i^2=1$, so that $\psi$ has degree 4. Clearly,
$\lambda=0$,
and $\deg\delta=3$, since the universal family contains exactly three singular
fibers.
The result follows.

With the next example in mind, we now make a general remark. Fix a sequence
$m_*=(0,m_1,\dots)$ and a positive integer $n$, and denote by $v$, $l$, and $g$
the
number of vertices, of edges, and the genus of any graph belonging to
$\mcomb{m_*}{n}$. The real dimension of $\wmn{m_*}{n}$ equals $l-n$, while $g$
is
given by $2-2g=v-l+n$.  It follows that $v\leq 0$, and hence $\wmn{m_*}{n}$ is
empty,
as soon as the real codimension of $\wmn{m_*}{n}$ equals or exceeds $4g-4+2n$.
In
particular, in this range, the formulas we are after would amount to expressing
certain
polynomial in the Mumford classes as linear combinations of boundary classes.
Moreover, if indeed these formulas were given by the DFIZ theorem, one could
conclude, by induction on the level, that all monomials in the Mumford classes
vanish
on
$\mgn{g}{n}$ in real codimension at least $4g-4+2n$. This is indeed true, as
follows from the observation by Harer (cf. [5], for instance) that $\mgn{g}{n}$
has the homotopy type of a CW-complex of dimension $4g-4+n$.

We next look at the codimension two classes $\wmn{m_*}{n}$ on $\mgnbar{1}{2}$.
The
remark we just made implies in particular that these are both zero. The
corresponding
conjectural formulas coming from the DFIZ theorem are
$$
0=120\kappa_2-6{\xi_{irr}}_*(\psi_1)-6
{\xi_{1,\emptyset}}_*(\psi_1\times 1)
+30{\xi_{1,\emptyset}}_*(\psi_1\times 1)\,,\leqno{(5.1)}
$$
$$
\eqalign{
0&=72\kappa_1^2-348\kappa_2
-6{\xi_{irr}}_*(\kappa_1)-12{\xi_{1,\emptyset}}_*(\kappa_1\times1)
+12{\xi_{irr}}_*(\psi_4)+12{\xi_{1,\emptyset}}_*(\psi_1\times 1)\cr
&\quad+{1\over
2}\sum_{p,P}{\xi_{(B,p,P)}}_*(1)+{1\over 4}\sum_{p,P}{\xi_{(D,p,P)}}_*(1)
-36{\xi_{1,\emptyset}}_*(\psi_1\times 1)\,.}\leqno{(5.2)}
$$
This last formula is a special case of formula (3.15). We have used the fact
that
graphs of type $A$ and $B$ cannot occur in our situation. It should also be
observed that the term corresponding to graph $B$ consists of a single summand,
for the two marked points must of necessity be on the rational tail. On the
other
hand, two distinct summands appear in the term corresponding to graph $D$. In
fact, this corresponds geometrically to two smooth rational curves joined at
two
points, each component carrying a marked point, and these can be labelled in
two
different ways. Now, using the computations of the two preceding examples, we
have that
$$
\eqalign{
&\int_{\mgnbar{1}{2}}{\xi_{irr}}_*(\kappa_1)=\int_{\mgnbar{0}{4}}\kappa_1=1
\,,\qquad
\int_{\mgnbar{1}{2}}{\xi_{irr}}_*(\psi_4)=\int_{\mgnbar{0}{4}}\psi_4=1\,,\cr
&\int_{\mgnbar{1}{2}}{\xi_{1,\emptyset}}_*(\psi_1\times
1)=\int_{\mgnbar{1}{1}}\psi_1=\frac{1}{24} \,,\qquad
\int_{\mgnbar{1}{2}}{\xi_{1,\emptyset}}_*(\kappa_1\times
1)=\int_{\mgnbar{1}{1}}\kappa_1=\frac{1}{24}\,,\cr
&\int_{\mgnbar{1}{2}}{\xi_{(B,p,P)}}_*(1)=
\int_{\mgnbar{0}{3}\times\mgnbar{0}{3}}1=1
\,,\qquad
\int_{\mgnbar{1}{2}}{\xi_{(D,p,P)}}_*(1)=
\int_{\mgnbar{0}{3}\times\mgnbar{0}{3}}1=1\,.
}
$$
On the other hand one has
$$
\int_{\mgnbar{1}{2}}\kappa_2=\frac{1}{24}\,,\qquad
\int_{\mgnbar{1}{2}}\kappa_1^2=\frac{1}{8}\,.
$$
This follows either from a simple algebro-geometric calculation or,
alternatively,
by noticing that the integrals to be computed are just
$$
\ang{\tau_0^2\tau_3}=\ang{\tau_1}=\frac{1}{24}
$$
and
$$
\ang{\tau_0^2\tau_2^2}-\ang{\tau_0^2\tau_3}=2\ang{\tau_0\tau_1\tau_2}-
\ang{\tau_1}=2\ang{\tau_0\tau_2}+2\ang{\tau_1\tau_1}-\ang{\tau_1}=
3\ang{\tau_1}=\frac{1}{8}\,.
$$
Substituting these values in the right-hand sides of (5.1) and (5.2) gives
zero, as
desired.

\vskip0.5truecm
We end this section with a few remarks. The first one concerns possible
generalizations of Lemma 2, and hence of Proposition 1, of section 4 to higher
codimension. It is clear that an essential ingredient in the proof of that
lemma is the
possibility of writing every degree two cohomology class as a linear
combination of standard ones. The analogue of this is only known to hold in
degree four, although it is a standard conjecture that it should in fact hold
in
every degree, provided the genus is sufficiently large. However, it would not
be without interest, and perhaps provable with the same methods we have
used in this section, that an analogue of Lemma 2 holds in all degrees,
provided attention is restricted only to those cohomology classes which can be
expressed as linear combinations of standard ones.

The second remark has to do with relations among standard classes. There is a
set of conjectures, due to Faber [unpublished], dealing with the relations that
the classes $\kappa_i$ satisfy in the rational cohomology of $\mgn{g}{n}$. It
is
known [9][10] that there are no such relations in degree less
than
$g/6$. For higher degrees Faber provides an explicit algebro-geometric recipe
to generate relations which, conjecturally, should yield all relations. It
occurred
to us that perhaps an alternative way of obtaining relations could be via a
recent result of Mulase [11] which states that the function $Z(t_*,s_*)$
satisfies the KdV hierarchy as a function of
$s_*$, for any fixed $t_*$. Making these equations explicit would yield
relations among the derivatives of $F$ with respect to the $s$ variables, and
we have explained how these could be translated into relations among the
$\kappa_i$.

Finally, it is clear that one needs to understand better the DFIZ theorem. In
particular, one should try and systematically explain the marvellous
cancellations that experimentally occur in all the cases we have been able to
compute.

\vskip1.5 truecm
\specialhead Appendix
\endspecialhead
Below are listed the expressions of the derivatives of $F$ with respect to the
$s$
variables in terms of derivatives with respect to $t$ variables that are
relevant to the
problem of expressing classes $\wmn{m_*}{n}$ in terms of algebro-geometric
classes, up
to weight 15. We recall that the weight of a partial derivative
$\prod(\partial/\partial
s_i)^{m_i}$ is defined to be $\sum m_i(2i+1)$. Of course the equalities below
hold {\it
only at} $s_*=\hat{s}_*=(0,1,0,\dots)$.
$$\eqalign{{\partial F\over\partial s_2}&=12\,{\partial F\over\partial t_2 } -
  {1\over 2}\,{\partial^2F\over\partial t_0^2} -
  {1\over 2}\,\left({\partial F\over\partial t_0 }\right)^2 }$$

$$\eqalign{{\partial F\over\partial s_3}&=120\,{\partial F\over\partial t_3 } -
  6\,{\partial^2F\over\partial t_0 \partial t_1 } -
  6\,{\partial F\over\partial t_1 }\,
   {\partial F\over\partial t_0 } +
  {5\over 4}\,{\partial F\over\partial t_0 } }$$
\vfill\break
$$\eqalign{{\partial F\over\partial s_4}&=1680\,{\partial F\over\partial t_4 }
-
  18\,{\partial^2F\over\partial t_1^2} -
  18\,{\left({\partial F\over\partial t_1 }\right)^2} -
  60\,{\partial^2F\over\partial t_0 \partial t_2 }  -
  60\,{\partial F\over\partial t_2 }\,
   {\partial F\over\partial t_0 } +
  {7\over 6}\,{\partial^3F\over\partial t_0^3}\cr &\quad +
  {7\over 2}\,{\partial F\over\partial t_0 }\,
      {\partial^2F\over\partial t_0^2}  +
  {7\over 6}\,\left({\partial F\over\partial t_0 }\right)^3 +
  {49\over 2}\,{\partial F\over\partial t_1 } -
  {{35}\over {96}} }$$

$$\eqalign{{\partial^2 F\over\partial s_2^2}&=144\,{\partial^2F\over\partial
t_2^2}
  - 840\,{\partial F\over\partial t_3 } -
  12\,{\partial^3F\over\partial t_0^2\partial t_2 }-
  24\,{\partial F\over\partial t_0 }\,
   {\partial^2F\over\partial t_0 \partial t_2 }  +
  24\,{\partial^2F\over\partial t_0 \partial t_1 } +
  24\,{\partial F\over\partial t_1 }\,
   {\partial F\over\partial t_0 }\cr &\quad +
  {1\over 4}\,{\partial^4F\over\partial t_0^4} +
  {\partial F\over\partial t_0 }\,{\partial^3F\over\partial t_0^3} +
  {1\over 2}\,\left({\partial^2F\over\partial t_0^2}\right)^2 +
  \left({\partial F\over\partial t_0 }\right)^2\,
      {\partial^2F\over\partial t_0^2} -
  3\,{\partial F\over\partial t_0 } }$$

$$\eqalign{{\partial F\over\partial s_5}&=30240\,{\partial F\over\partial t_5 }
-
  360\,{\partial^2F\over\partial t_1 \partial t_2 } -
  360\,{\partial F\over\partial t_2 }\,{\partial F\over\partial t_1 } -
  840\,{\partial^2F\over\partial t_0 \partial t_3 } -
  840\,{\partial F\over\partial t_3 }\,
   {\partial F\over\partial t_0 }\cr &\quad +
  27\,{\partial^3F\over\partial t_0^2\partial t_1 } +
  27\,{\partial F\over\partial t_1 }\,{\partial^2F\over\partial t_0^2} +
  54\,{\partial F\over\partial t_0 }\,
   {\partial^2F\over\partial t_0 \partial t_1 } +
  27\,{\partial F\over\partial t_1 }\,{\left({\partial F\over\partial t_0
}\right)^2} +
  585\,{\partial F\over\partial t_2 }\cr &\quad -
  {105\over 8}\,{\partial^2F\over\partial t_0^2} -
  {105\over 8}\,{\left({\partial F\over\partial t_0 }\right)^2} }$$

$$\eqalign{{\partial^2 F\over\partial s_2\partial s_3}&=
   1440\,{\partial^2F\over\partial t_2 \partial t_3 } -
   15120\,{\partial F\over\partial t_4 } -
   60\,{\partial^3F\over\partial t_0^2\partial t_3 } -
  120\,{\partial F\over\partial t_0 }\,{\partial^2F\over\partial t_0 \partial
t_3 }  -
  72\,{\partial^3F\over\partial t_0 \partial t_1 \partial t_2 }\cr &\quad -
  72\,{\partial F\over\partial t_1 }\,
   {\partial^2F\over\partial t_0 \partial t_2 } -
  72\,{\partial^2F\over\partial t_1 \partial t_2 }\,
   {\partial F\over\partial t_0 }  +
  90\,{\partial^2F\over\partial t_1^2} +
  90\,{\left({\partial F\over\partial t_1 }\right)^2} +
  375\,{\partial^2F\over\partial t_0 \partial t_2 }\cr &\quad +
  360\,{\partial F\over\partial t_2 }\,
   {\partial F\over\partial t_0 }  +
  3\,{\partial^4F\over\partial t_0^3\partial t_1 } +
  6\,{\partial^2F\over\partial t_0 \partial t_1 }\,
   {\partial^2F\over\partial t_0^2} +
  9\,{\partial F\over\partial t_0 }\,{\partial^3F\over\partial t_0^2\partial
t_1 }  +
  3\,{\partial F\over\partial t_1 }\,{\partial^3F\over\partial t_0^3}\cr &\quad
+
  6\,{\partial F\over\partial t_1 }\,{\partial F\over\partial t_0 }\,
   {\partial^2F\over\partial t_0^2} +
  6\,{\left({\partial F\over\partial t_0}\right)^2}\,
   {\partial^2F\over\partial t_0 \partial t_1 } -
  {45\over 8}\,{\partial^3F\over\partial t_0^3} -
  {65\over 4}\,{\partial F\over\partial t_0 }\,{\partial^2F\over\partial t_0^2}
 -
  5\,{\left({\partial F\over\partial t_0 }\right)^3}\cr &\quad -
  {165\over 2}\,{\partial F\over\partial t_1 } +
  {{29}\over {32}} }$$

$$\eqalign{{\partial F\over\partial s_6}&=
  665280\,{\partial F\over\partial t_6 } -
  1800\,{\partial^2F\over\partial t_2^2} -
  1800\,{\left({\partial F\over\partial t_2}\right)^2} -
  5040\,{\partial^2F\over\partial t_1 \partial t_3 } -
  5040\,{\partial F\over\partial t_3 }\,
   {\partial F\over\partial t_1 }\cr &\quad -
  15120\,{\partial^2F\over\partial t_0 \partial t_4 } -
  15120\,{\partial F\over\partial t_4 }\,{\partial F\over\partial t_0 } +
  16170\,{\partial F\over\partial t_3 } +
  198\,{\partial^3F\over\partial t_0 \partial t_1^2} +
  396\,{\partial F\over\partial t_1 }\,{\partial^2F\over\partial t_0 \partial
t_1 }\cr
&\quad +
  198\,{\partial^2F\over\partial t_1^2}\,{\partial F\over\partial t_0 } +
  198\,{\left({\partial F\over\partial t_1 }\right)^2}\,{\partial
F\over\partial t_0 } +
  330\,{\partial^3F\over\partial t_0^2\partial t_2 } +
  330\,{\partial F\over\partial t_2 }\,
   {\partial^2F\over\partial t_0^2}\cr &\quad +
  660\,{\partial F\over\partial t_0 }\,
   {\partial^2F\over\partial t_0 \partial t_2 } +
  330\,{\partial F\over\partial t_2 }\,{\left({\partial F\over\partial t_0
}\right)^2} -
  {33\over 8}\,{\partial^4F\over\partial t_0^4} -
  {33\over 2}\,{\partial F\over\partial t_0 }\,
      {\partial^3F\over\partial t_0^3} -
  {99\over 8}\,{\left({\partial^2F\over\partial t_0^2}\right)^2}\cr &\quad -
  {99\over 4}\,{\left({\partial F\over\partial t_0 }\right)^2}\,
      {\partial^2F\over\partial t_0^2} -
  {33\over 8}\,{\left({\partial F\over\partial t_0 }\right)^4} -
  {891\over 2}\,{\partial^2F\over\partial t_0 \partial t_1 } -
  {891\over 2}\,{\partial F\over\partial t_1 }\,{\partial F\over\partial t_0 }
+
  {1155\over 32}\,{\partial F\over\partial t_0 } }$$

$$\eqalign{{\partial^2 F\over\partial s_2\partial s_4}&=
  20160\,{\partial^2F\over\partial t_2 \partial t_4 } -
  332640\,{\partial F\over\partial t_5 } -
  216\,{\partial^3F\over\partial t_1^2\partial t_2 } -
  432\,{\partial F\over\partial t_1 }\,{\partial^2F\over\partial t_1 \partial
t_2 } -
  840\,{\partial^3F\over\partial t_0^2\partial t_4 }\cr &\quad -
  1680\,{\partial F\over\partial t_0 }\,{\partial^2F\over\partial t_0 \partial
t_4 } -
  720\,{\partial^3F\over\partial t_0 \partial t_2^2} -
  720\,{\partial^2F\over\partial t_2^2}\,
   {\partial F\over\partial t_0 } -
  720\,{\partial F\over\partial t_2 }\,
   {\partial^2F\over\partial t_0 \partial t_2 }\cr &\quad +
  6720\,{\partial^2F\over\partial t_0 \partial t_3 } +
  6720\,{\partial F\over\partial t_3 }\,
   {\partial F\over\partial t_0 } +
  2814\,{\partial^2F\over\partial t_1 \partial t_2 } +
  2520\,{\partial F\over\partial t_2 }\,
   {\partial F\over\partial t_1 }\cr &\quad +
  9\,{\partial^4F\over\partial t_0^2\partial t_1^2} +
  18\,{\partial F\over\partial t_1 }\,{\partial^3F\over\partial t_0^2\partial
t_1 } +
  18\,{\left({\partial^2F\over\partial t_0 \partial t_1 }\right)^2} +
  18\,{\partial F\over\partial t_0 }\,{\partial^3F\over\partial t_0 \partial
t_1^2}\cr
&\quad +
  36\,{\partial F\over\partial t_1 }\,{\partial F\over\partial t_0 }\,
   {\partial^2F\over\partial t_0 \partial t_1 } +
  44\,{\partial^4F\over\partial t_0^3\partial t_2 } +
  30\,{\partial F\over\partial t_2 }\,
   {\partial^3F\over\partial t_0^3} +
  102\,{\partial^2F\over\partial t_0 \partial t_2 }\,{\partial^2F\over\partial
t_0^2}
  \cr &\quad +
  132\,{\partial F\over\partial t_0 }\,
   {\partial^3F\over\partial t_0^2\partial t_2 } +
  102\,{\left({\partial F\over\partial t_0 }\right)^2}\,
   {\partial^2F\over\partial t_0 \partial t_2 } +
  60\,{\partial F\over\partial t_2 }\,{\partial F\over\partial t_0 }\,
   {\partial^2F\over\partial t_0^2} -
   2835\,{\partial F\over\partial t_2 }\cr &\quad -
  {637\over 4}\,{\partial^3F\over\partial t_0^2\partial t_1 } -
  147\,{\partial F\over\partial t_1 }\,{\left({\partial F\over\partial t_0
}\right)^2} -
  147\,{\partial F\over\partial t_1 }\,{\partial^2F\over\partial t_0^2} -
  {637\over 2}\,{\partial F\over\partial t_0 }
  \,{\partial^2F\over\partial t_0 \partial t_1 }\cr &\quad -
  {7\over 12}\,{\partial^5F\over\partial t_0^5} -
  {35\over 12}\,{\partial F\over\partial t_0 }\,
      {\partial^4F\over\partial t_0^4} -
  {21\over 4}\,{\partial^2F\over\partial t_0^2}\,{\partial^3F\over\partial
t_0^3} -
  7\,{\partial F\over\partial t_0 }\,
      {\left({\partial^2F\over\partial t_0^2}\right)^2}\cr &\quad -
  {21\over 4}\,{\left({\partial F\over\partial t_0 }\right)^2}\,
      {\partial^3F\over\partial t_0^3} -
  {7\over 2}\,{\left({\partial F\over\partial t_0 }\right)^3}\,
      {\partial^2F\over\partial t_0^2} +
  {385\over 8}\,{\partial^2F\over\partial t_0^2} +
  {385\over 8}\,{\left({\partial F\over\partial t_0 }\right)^2} }$$

$$\eqalign{{\partial^2 F\over\partial s_3^2}&=
  14400\,{\partial^2F\over\partial t_3^2} -
  332640\,{\partial F\over\partial t_5 } -
  1440\,{\partial^3F\over\partial t_0 \partial t_1 \partial t_3 } -
  1440\,{\partial^2F\over\partial t_1 \partial t_3 }\,
   {\partial F\over\partial t_0 } -
  1440\,{\partial F\over\partial t_1 }\,
   {\partial^2F\over\partial t_0 \partial t_3 }\cr &\quad +
  7020\,{\partial^2F\over\partial t_0 \partial t_3 } +
  6720\,{\partial F\over\partial t_3 }\,
   {\partial F\over\partial t_0 } +
  2160\,{\partial^2F\over\partial t_1 \partial t_2 } +
  2160\,{\partial F\over\partial t_2 }\,
   {\partial F\over\partial t_1 } +
  36\,{\partial^4F\over\partial t_0^2\partial t_1^2}\cr &\quad +
  36\,{\left({\partial^2F\over\partial t_0 \partial t_1 }\right)^2} +
  36\,{\partial^2F\over\partial t_1^2}\,
   {\partial^2F\over\partial t_0^2} +
  72\,{\partial F\over\partial t_1 }\,
   {\partial^3F\over\partial t_0^2\partial t_1 } +
  72\,{\partial F\over\partial t_0 }\,
   {\partial^3F\over\partial t_0 \partial t_1^2}\cr &\quad +
  36\,{\left({\partial F\over\partial t_1 }\right)^2}\,
   {\partial^2F\over\partial t_0^2} +
  72\,{\partial F\over\partial t_1 }\,
   {\partial F\over\partial t_0 }\,
   {\partial^2F\over\partial t_0 \partial t_1 } +
  36\,{\partial^2F\over\partial t_1^2}\,
   {\left({\partial F\over\partial t_0 }\right)^2} -
  2400\,{\partial F\over\partial t_2 }\cr &\quad -
  165\,{\partial^3F\over\partial t_0^2\partial t_1 } -
  165\,{\partial F\over\partial t_1 }\,
   {\partial^2F\over\partial t_0^2} -
  150\,{\partial F\over\partial t_1 }\,
   {\left({\partial F\over\partial t_0 }\right)^2} -
  315\,{\partial F\over\partial t_0 }\,
   {\partial^2F\over\partial t_0 \partial t_1 }\cr &\quad +
  {725\over 16}\,{\partial^2F\over\partial t_0^2} +
  {175\over 4}\,{\left({\partial F\over\partial t_0 }\right)^2} }$$

$$\eqalign{{\partial F\over\partial s_7}&=17297280\,{\partial F\over\partial
t_7 }
-
  50400\,{\partial^2F\over\partial t_2 \partial t_3 } -
  50400\,{\partial F\over\partial t_3 }\,
   {\partial F\over\partial t_2 } -
  90720\,{\partial^2F\over\partial t_1 \partial t_4 } -
  90720\,{\partial F\over\partial t_4 }\,
   {\partial F\over\partial t_1 }\cr &\quad -
  332640\,{\partial^2F\over\partial t_0 \partial t_5 } -
  332640\,{\partial F\over\partial t_5 }\,
   {\partial F\over\partial t_0 } +
  468\,{\partial^3F\over\partial t_1^3} +
  1404\,{\partial F\over\partial t_1 }\,
   {\partial^2F\over\partial t_1^2} +
  468\,{\left({\partial F\over\partial t_1 }\right)^3}\cr &\quad +
  4680\,{\partial^3F\over\partial t_0 \partial t_1 \partial t_2 } +
  4680\,{\partial F\over\partial t_1 }\,
   {\partial^2F\over\partial t_0 \partial t_2 } +
  4680\,{\partial F\over\partial t_2 }\,
   {\partial^2F\over\partial t_0 \partial t_1 } +
  4680\,{\partial^2F\over\partial t_1 \partial t_2 }\,
   {\partial F\over\partial t_0 }\cr &\quad +
  4680\,{\partial F\over\partial t_2 }\,
   {\partial F\over\partial t_1 }\,
   {\partial F\over\partial t_0 } +
  5460\,{\partial^3F\over\partial t_0^2\partial t_3 } +
  5460\,{\partial F\over\partial t_3 }\,
   {\partial^2F\over\partial t_0^2} +
  10920\,{\partial F\over\partial t_0 }\,
   {\partial^2F\over\partial t_0 \partial t_3 }\cr &\quad +
  5460\,{\partial F\over\partial t_3 }\,
   {\left({\partial F\over\partial t_0 }\right)^2} +
  507780\,{\partial F\over\partial t_4 } -
  10725\,{\partial^2F\over\partial t_0 \partial t_2 } -
  10725\,{\partial F\over\partial t_2 }\,
   {\partial F\over\partial t_0 }\cr &\quad -
  {5577\over 2}\,{\partial^2F\over\partial t_1^2} -
  {5577\over 2}\,{\left({\partial F\over\partial t_1 }\right)^2} -
  143\,{\partial^4F\over\partial t_0^3\partial t_1 } -
  143\,{\partial F\over\partial t_1 }\,
   {\partial^3F\over\partial t_0^3} -
  429\,{\partial^2F\over\partial t_0 \partial t_1 }\,
   {\partial^2F\over\partial t_0^2}\cr &\quad -
  429\,{\partial F\over\partial t_0 }\,
   {\partial^3F\over\partial t_0^2\partial t_1 } -
  429\,{\left({\partial F\over\partial t_0 }\right)^2}\,
   {\partial^2F\over\partial t_0 \partial t_1 } -
  429\,{\partial F\over\partial t_1 }\,
   {\partial F\over\partial t_0 }\,
   {\partial^2F\over\partial t_0^2} -
  143\,{\partial F\over\partial t_1 }\,
   {\left({\partial F\over\partial t_0 }\right)^3}\cr &\quad +
  {1001\over 8}\,{\partial^3F\over\partial t_0^3} +
  {3003\over 8}\,{\partial F\over\partial t_0 }\,
   {\partial^2F\over\partial t_0^2} +
  {1001\over 8}\,{\left({\partial F\over\partial t_0 }\right)^3} +
  {27027\over 16}\,{\partial F\over\partial t_1 } -
  {{5005}\over {384}} }$$

$$\eqalign{{\partial^3 F\over\partial s_2^3}&=1728\,{\partial^3F\over\partial
t_2^3} -
  30240\,{\partial^2F\over\partial t_2 \partial t_3 } -
  216\,{\partial^4F\over\partial t_0^2\partial t_2^2} -
  432\,{\left({\partial^2F\over\partial t_0 \partial t_2 }\right)^2} -
  432\,{\partial F\over\partial t_0 }\,{\partial^3F\over\partial t_0 \partial
t_2^2}\cr
&\quad +
  864\,{\partial^3F\over\partial t_0 \partial t_1 \partial t_2 } +
  864\,{\partial F\over\partial t_1 }\,
   {\partial^2F\over\partial t_0 \partial t_2 } +
  864\,{\partial^2F\over\partial t_1 \partial t_2 }\,
   {\partial F\over\partial t_0 } +
  1260\,{\partial^3F\over\partial t_0^2\partial t_3 }\cr &\quad +
  2520\,{\partial F\over\partial t_0 }\,{\partial^2F\over\partial t_0 \partial
t_3 } +
  9\,{\partial^5F\over\partial t_0^4\partial t_2 } +
  36\,{\partial^2F\over\partial t_0 \partial t_2 }\,{\partial^3F\over\partial
t_0^3} +
  36\,{\partial F\over\partial t_0 }\,{\partial^4F\over\partial t_0^3\partial
t_2 } +
  36\,{\partial^2F\over\partial t_0^2}\,{\partial^3F\over\partial t_0^2\partial
t_2}
\cr &\quad +
  72\,{\partial F\over\partial t_0 }\,{\partial^2F\over\partial t_0 \partial
t_2 }\,
   {\partial^2F\over\partial t_0^2} +
  36\,{\left({\partial F\over\partial t_0 }\right)^2}
  \,{\partial^3F\over\partial t_0^2\partial t_2 } +
  151200\,{\partial F\over\partial t_4 } -
  576\,{\partial^2F\over\partial t_1^2} -
  576\,{\left({\partial F\over\partial t_1 }\right)^2}\cr &\quad
 -
  2628\,{\partial^2F\over\partial t_0 \partial t_2 } -
  2520\,{\partial F\over\partial t_2 }\,
   {\partial F\over\partial t_0 } -
  36\,{\partial^4F\over\partial t_0^3\partial t_1 } -
  108\,{\partial F\over\partial t_0}\,
   {\partial^3F\over\partial t_0^2\partial t_1 } -
  36\,{\partial F\over\partial t_1 }\,{\partial^3F\over\partial t_0^3}\cr
&\quad -
  72\,{\partial^2F\over\partial t_0 \partial t_1 }\,
   {\partial^2F\over\partial t_0^2} -
  72\,{\partial F\over\partial t_1 }\,
   {\partial F\over\partial t_0 }\,{\partial^2F\over\partial t_0^2} -
  72\,{\left({\partial F\over\partial t_0 }\right)^2}
  \,{\partial^2F\over\partial t_0 \partial t_1 } -
  {1\over 8}\,{\partial^6F\over\partial t_0^6} -
  {3\over 4}\,{\partial F\over\partial t_0 }\,
      {\partial^5F\over\partial t_0^5}\cr &\quad -
  {3\over 2}\,{\partial^2F\over\partial t_0^2}\,{\partial^4F\over\partial
t_0^4} -
  {5\over 4}\,{\left({\partial^3F\over\partial t_0^3}\right)^2} -
  6\,{\partial F\over\partial t_0 }\,{\partial^2F\over\partial t_0^2}\,
   {\partial^3F\over\partial t_0^3} -
  {\left({\partial^2F\over\partial t_0^2}\right)^3} -
  {3\over 2}\,{\left({\partial F\over\partial t_0 }\right)^2}\,
      {\partial^4F\over\partial t_0^4}\cr &\quad -
  3\,{\left({\partial F\over\partial t_0 }\right)^2}\,
   {\left({\partial^2F\over\partial t_0^2}\right)^2} -
  {\left({\partial F\over\partial t_0 }\right)^3}\,{\partial^3F\over\partial
t_0^3} +
  {63\over 2}\,{\partial^3F\over\partial t_0^3} +
  90\,{\partial F\over\partial t_0 }\,
   {\partial^2F\over\partial t_0^2} +
  27\,{\left({\partial F\over\partial t_0 }\right)^3}\cr &\quad +
  378\,{\partial F\over\partial t_1 } -
  {{63}\over {20}} }$$

\vfill\break

\enddocument